\documentclass[11pt]{article}
\usepackage[T1]{fontenc}
\usepackage[margin=2.5cm]{geometry}
\usepackage{graphicx}
\usepackage[utf8]{inputenc}
\usepackage{titletoc}

\usepackage{suffix}
\usepackage[svgnames,x11names]{xcolor}
\usepackage{xfrac,xifthen,xparse,xpunctuate,xspace}

\usepackage{amsmath,amsfonts,amssymb,dsfont,latexsym,mathtools,mathrsfs}

\DeclareMathOperator*{\argmax}{arg\,max}
\DeclareMathAlphabet\mathpzc{OT1}{pzc}{m}{it}

\usepackage{accents}

\usepackage{tikz}
\usepackage{caption,subcaption}
\usetikzlibrary{calc}

\DeclareCaptionLabelFormat{figcap}{#1~#2}
\DeclareCaptionLabelFormat{subfigcap}{Figure~\thefigure#2}
\usepackage[amsmath,hyperref,thmmarks]{ntheorem}
\usepackage{amssymb,dsfont,enumitem,mathrsfs,thmtools,upgreek}

\newcommand\proofSymbol{\ensuremath{_\blacksquare}}
\newcommand\qedSymbol{\proofSymbol}
\newcommand\qedhere{\ifmmode\qed\else\hfill\proofSymbol\fi}

\makeatletter
\def\theorem@checkbold{}
\newtheoremstyle{stmstyle}%
  {\item{\theorem@headerfont ##1\ ##2\theorem@separator}\enspace}%
  {\item{\theorem@headerfont ##1\ ##2\ {\normalfont(##3)}%
    \theorem@separator}\enspace}
\newtheoremstyle{prfstyle}%
  {\item{\theorem@headerfont ##1\theorem@separator}\enspace}%
  {\item{\theorem@headerfont ##3\theorem@separator}\enspace}
\makeatother
\theoremseparator{.}
\theorembodyfont{\upshape}

\theoremheaderfont{\itshape}
\theoremsymbol{\qedSymbol}
\declaretheorem[name=Proof,style=prfstyle]{prfenv}

\theoremheaderfont{\sffamily\bfseries}
\theoremsymbol{}

\declaretheorem[name=Proposition,style=stmstyle,numberlike=lemenv]%
  {proenv}

\usepackage{tcolorbox}
\tcbuselibrary{theorems}

\newtcbtheorem[number within=section]{ctheorem}{Theorem}%
{colback=ared!5,colframe=ared!35!black,fonttitle={\bfseries\sffamily}}{thm}

\usepackage[
  colorlinks=true,
  linkcolor=MediumBlue,
  citecolor=Firebrick2,
  hypertexnames=false
]{hyperref}
\usepackage[capitalise,compress,nameinlink]{cleveref}

\crefname{lemenv}{Lemma}{Lemmas}
\crefname{proenv}{Proposition}{Propositions}
\crefname{thmenv}{Theorem}{Theorems}
\crefname{tcb@cnt@ctheorem}{Theorem}{Theorems}
\crefname{corenv}{Corollary}{Corollaries}
\crefname{remenv}{Remark}{Remarks}
\crefname{exmenv}{Example}{Examples}
\crefname{figure}{Figure}{Figures}
\crefname{dfnenv}{Definition}{Definitions}

\newlist{assump}{enumerate}{1}
\setlist[assump]{label=\sffamily(H\arabic*)}
\crefname{assumpi}{Assumption}{Assumptions}

\edef\envlabel{}
\newcommand\shortenv[2][]{%
\def\@tempa{#1}%
\DeclareDocumentEnvironment{\ifx\@tempa\@empty\relax#1\else #2\fi}%
{ o g }{\ifvmode\else\unskip\fi%
  \IfValueTF{##1}{\begin{#1env}[##1]}{\begin{#1env}}%
  \IfValueT{##2}{\edef\envlabel{##2}\label{\envlabel}}}%
{\end{#1env}}%
}
\shortenv[lem]{lemma}
\shortenv[pro]{proposition}
\shortenv[thm]{theorem}
\shortenv[cor]{corollary}
\shortenv[rem]{remark}
\shortenv[exm]{example}
\shortenv[dfn]{definition}

\DeclareDocumentEnvironment{proof}{ o g }{%
  \IfValueTF{#1}{\begin{prfenv}[#1]}{\begin{prfenv}}}
{\end{prfenv}}

\makeatletter
\let\c@tcb@cnt@ctheorem\c@lemenv
\newcommand{\leqnomode}{\tagsleft@true}
\newcommand{\reqnomode}{\tagsleft@false}

\def\@MRExtract#1 #2!{#1}     
\newcommand{\MR}[1]{
  \xdef\@MRSTRIP{\@MRExtract#1 !}%
  \href{http://www.ams.org/mathscinet-getitem?mr=\@MRSTRIP}{MR-\@MRSTRIP}}
\newcommand{\ARXIV}[1]{\href{http://arXiv.org/abs/#1}{arXiv:#1}}
\makeatother

\catcode`@=11
\def\nvphantom{\v@true\h@false\nph@nt}
\def\nhphantom{\v@false\h@true\nph@nt}
\def\nphantom{\v@true\h@true\nph@nt}
\def\nph@nt{\ifmmode\def\next{\mathpalette\nmathph@nt}%
  \else\let\next\nmakeph@nt\fi\next}
\def\nmakeph@nt#1{\setbox\z@\hbox{#1}\nfinph@nt}
\def\nmathph@nt#1#2{\setbox\z@\hbox{$\m@th#1{#2}$}\nfinph@nt}
\def\nfinph@nt{\setbox\tw@\null
  \ifv@ \ht\tw@\ht\z@ \dp\tw@\dp\z@\fi
  \ifh@ \wd\tw@-\wd\z@\fi \box\tw@}

\newcommand\NN{\mathbb N}
\newcommand\ZZ{\mathbb Z}
\newcommand\RR{\mathbb R}
\newcommand\EE{\mathbb E}
\newcommand\PP{\mathbb P}
\newcommand\II{\mathds 1}
\newcommand\eps{\varepsilon}
\newcommand\es{\varnothing}

\newcommand\Leb{\mathrm L}
\newcommand\Co{\mathcal C}

\newcommand\defeq{\coloneqq}
\newcommand\eqdef{\eqqcolon}
\newcommand\equiveq{\mathrel{\rlap{\raisebox{0.3ex}{$\cdot$}}%
  \raisebox{-0.3ex}{$\cdot$}}\equiv}
\newcommand\eqlaw{\overset d=}

\newcommand\dd{\mathrm d}
\newcommand\ubar[1]{\text{\b{$#1$}}}
\newcommand\bs[1]{\boldsymbol{\mathbf{#1}}}

\newcommand\oo[2]{\mathopen{(}#1,#2\mathclose{)}}
\newcommand\oc[2]{\mathopen{(}#1,#2\mathclose{]}}
\newcommand\co[2]{\mathopen{[}#1,#2\mathclose{)}}
\newcommand\cc[2]{\mathopen{[}#1,#2\mathclose{]}}

\newcommand\MP{\mathscr P}
\newcommand\DSS[1][2]{\ell^{#1\downarrow}}
\newcommand\BRW{\mathcal Z}
\newcommand\BRWb{\widetilde{\mathcal Z}}
\newcommand\cf{\kappa}
\newcommand\cfb{\widetilde\kappa}
\newcommand\cffr{{\omega_-}}
\newcommand\cfsr{{\omega_+}}
\newcommand\cem{\partial}
\newcommand\HUT{\mathbb U}
\newcommand\BT{\mathbb B}
\newcommand\bt{\BT(\eps)}
\newcommand\CS[1][X]{\mathcal #1}
\newcommand\RL{\mathcal L}
\newcommand\tc{{\xi^-}}
\newcommand\TC{\widehat{\mathcal X}}
\newcommand\Em[1][x]{\widehat{\mathcal{E}}^{-}_{#1}}
\newcommand\Pm[1][x]{\widehat{\mathcal{P}}^{-}_{#1}}
\newcommand\Ex[1][x]{\mathcal{E}_{#1}}
\newcommand\Px[1][x]{\mathcal{P}_{#1}}
\newcommand\Eb[1][x]{\mathbf{E}_{#1}}
\newcommand\Pb[1][x]{\mathbf{P}_{#1}}

\newcommand\fc[2][i]{{\mathpzc{x}\hspace{0.25pt}}_{#1,#2}}
\newcommand\ft[1][\RL]{\mathpzc{T}_{#1}(\eps)}
\newcommand\FM{\mathpzc{M}^-(\eps)}

\author{Benjamin Dadoun\thanks{Institut für Mathematik, 
  Universität Zürich, 
  Winterthurerstrasse 190, 
  CH-8057 Zürich, Switzerland.\hfill\eject
  Email: benjamin.dadoun@math.uzh.ch\vspace{.2em}\hfill\eject
\textit{{Keywords:}}  Growth-fragmentation,
self-similarity, additive martingale\hfill\eject
}}
\title{\scshape%
  Asymptotics of self-similar\\growth-fragmentation processes}
\date{}

\begin{document}
\maketitle
\begin{abstract}
Markovian growth-fragmentation processes introduced
  in~\cite{Bertoin16,Bertoin17} extend the pure-fragmentation model
  by allowing the fragments to grow larger or
  smaller between dislocation events. What becomes of the known
  asymptotic behaviors of self-similar pure
  fragmentations~\cite{Bertoin03,Gnedin04,Martinez05,Rouault05} when
  growth is added to the fragments is a natural question that we
  investigate in this paper. Our results involve the terminal value of
  some additive martingales whose uniform integrability is an
  essential requirement. Dwelling first on the homogeneous case~\cite{Bertoin16},
  we exploit the connection with branching random walks and in
  particular the martingale convergence of
  Biggins~\cite{Biggins77,Biggins92} to derive precise asymptotic
  estimates. The self-similar case~\cite{Bertoin17} is treated in a
  second part; under the so called Malthusian hypotheses and with
  the help of several martingale-flavored features recently developed
  in~\cite{Curien18}, we obtain limit theorems for empirical
  measures of the fragments.
\end{abstract}

\section{Introduction}
Fragmentation processes are meant to describe the evolution of an object
which is subject to random and repeated dislocations over time. The way
the mass is spread into smaller fragments during a dislocation event is
usually given by a (random) \emph{mass-partition}, that is an element of
the space
\begin{equation}
  \MP\defeq\left\{\bs p\defeq(p_i,\,i\in\NN)
  \colon p_1\ge p_2\ge\cdots\ge0\enspace\text{and}\enspace%
    \sum_{i=1}^{\infty}p_i\le1\right\}\!,%
  \label{eq:A_Newton}
\end{equation}
where the total mass need not be conserved, i.e.\ a positive proportion
$1-\sum_{i\ge1}p_i$ may disintegrate into
dust. The first probabilistic models of fragmentations go back at least
to Kolmogorov~\cite{Kolmogorov41}. Roughly, Kolmogorov imagined a
discrete branching system in\vadjust{\goodbreak} which particles get fragmented according
to a conservative distribution~$\nu$ on~$\MP$ and in a
\emph{homogeneous} manner, that is to say the rate at which a particle
splits does not depend on its mass. Under this essential assumption of
homogeneity, Kolmogorov showed that a simple rescaling of the empirical
measure of the logarithms of the fragments converges with probability
one toward the Gaussian distribution. Later, a student of his,
Filippov~\cite{Filippov61} investigated mass-dependent dislocation rates
and more precisely the \emph{self-similar} case, in the sense that a
particle with size~$m$ splits at speed~$m^\alpha$ for some fixed
constant $\alpha\in\RR$ (the homogeneous case then corresponds to
$\alpha=0$). Most notably he discovered a limit theorem for a weighted
version of the empirical measure of the fragments when $\alpha>0$.
The special but common \emph{binary} situation, where particles always
split into two smaller fragments, has been emphasized by Brennan and
Durrett~\cite{Brennan86,Brennan87}, and later reconsidered by
Baryshnikov and Gnedin~\cite{Gnedin01} in some variant of the car
packing problem. Further extensions and other asymptotic properties in
the non-conservative case have also been derived by Bertoin
and Gnedin~\cite{Gnedin04} by means of complex analysis and contour
integrals.

In the 2000s (see \cite[Chapters 1-3]{Bertoin06} for a comprehensive
summary), Bertoin extended and theorized the construction of general
fragmentation processes in continuous time. In particular the
\emph{dislocation measure}~$\nu$ need no longer be a probability
distribution, as there is only the integrability requirement
\begin{equation}
  \int_{\MP}(1-p_1)\,\nu(\dd\bs p)<\infty.%
  \label{eq:A_Archimedes}
\end{equation}
While permitting infinite dislocation rates (so infinitely
many dislocation events may occur in a bounded time interval), this
condition prevents the total mass from being immediately shattered into
dust and leads to a nondegenerate fragmentation process
$\bs X(t)\defeq(X_1(t),X_2(t),\ldots),\,t\ge0,$ with values in~$\MP$.
When $\alpha=0$, fragmentation processes can be related (via a simple
logarithmic transformation) to \emph{branching random walks}, for which
fruitful literature is available, see e.g.\ the works of Biggins and
Uchiyama~\cite{Biggins77,Uchiyama82,Biggins92}, and~\cite{Shi15}.
Especially, \emph{additive martingales}, which are processes of the form
\begin{equation}
  \EE\!\left[\sum_{i=1}^\infty X_i^q(t)\right]^{-1}
    \sum_{i=1}^\infty X_i^q(t),\qquad t\ge0,%
  \label{eq:A_Gauss}
\end{equation}
for some parameter $q>0$, play a key role and the question of their
uniform integrability inquired by Biggins has successfully led to the
asymptotic behavior of homogeneous conservative
fragmentations~\cite{Rouault04,Rouault05}. More generally, in the
self-similar case, some specific so called \emph{Malthusian hypotheses}
guarantee the existence of an intrinsic martingale associated with the
fragmentation and whose convergence again yields many interesting
asymptotic results. Among others the results of Kolmogorov and Filippov
have been revisited~\cite{Bertoin03}, applying known statistics of
self-similar Markov processes to the process of a randomly tagged
fragment.

More recently, Bertoin~\cite{Bertoin16,Bertoin17} introduced a new
type of fragmentation processes in which the fragments are allowed to
grow during their lifetimes. We expect that most of the aforementioned
asymptotic properties extend to these \emph{growth-fragmentation}
processes, and it is the main purpose of the present work to derive
some of them. We shall first give a bit more description and explain
why our task is not completely straightforward. Like in the pure (i.e.\
without growth) setting, we are interested in the process which
describes the (sizes of the) fragments as time passes. For homogeneous
growth-fragmentations, namely the \emph{compensated fragmentations}
of~\cite{Bertoin16}, the basic prototype is simply a dilated homogeneous
fragmentation, that is a pure homogeneous fragmentation affected by a
deterministic exponential drift. However, there exist much more general
compensated fragmentations, where the dislocation measure~$\nu$ has only
to fulfill
\begin{equation}
 \int_\MP(1-p_1)^2\,\nu(\dd\bs p)<\infty,\label{eq:A_Euler}
\end{equation}
so that the process is nondegenerate and can still be encoded at any
time by a non-increasing null sequence. Condition~\eqref{eq:A_Euler} is
weaker than the necessary and sufficient condition~\eqref{eq:A_Archimedes}
for~$\nu$ to be the dislocation measure of a homogeneous fragmentation,
and both are reminiscent of those concerning the jump intensities of
L\'{e}vy processes, respectively subordinators. Incidentally, it was the
main motivation of~\cite{Bertoin16} to establish that, just like the
L\'{e}vy--It\=o construction of L\'{e}vy processes in terms of compensated
Poisson integrals, compensated fragmentations naturally arise as limits
of suitably dilated homogeneous fragmentations
\cite[Theorem~2]{Bertoin16}. Though asymptotic properties of pure
homogeneous fragmentations immediately transfer to the dilated ones,
extending them to general compensated fragmentations would correspond to
interchanging two limits, which does not seem obvious at first sight.
This is without to mention the self-similar case, that is for the
growth-fragmentations in~\cite{Bertoin17}, where things look even more
complicated.

There, and unlike the compensated fragmentations which are constructed
directly as processes in time, the self-similar \emph{cell systems} are
rather built from a genealogical point of view: roughly, the (size of
the) mother cell evolves like a Markov process on the positive
half-line where each negative jump~$-y$ is interpreted as a splitting
event, giving birth to a daughter cell with initial size~$y$ and which
then grows independently of the mother particle and according to the
same dynamics, i.e.\ producing in turn granddaughters, and so on. Bertoin
focused in particular on the situation where the associated
growth-fragmentation process $\bs X\defeq(\bs X(t),\,t\ge0)$, that is the
process recording the sizes of all alive cells in the system, fulfills
a self-similarity property, namely when there exists $\alpha\in\RR$ such
that for each $x>0$, the process $(x\bs X(x^\alpha t),\,t\ge0)$ has the
same law as~$\bs X$ started from a cell whose initial size is~$x$. In
the homogeneous case $\alpha=0$, these growth-fragmentations correspond
to the compensated fragmentations of~\cite{Bertoin16} for which the
dislocation measure is binary, see
\cite[Proposition~3]{Bertoin17}. In the self-similar case $\alpha<0$,
they have been proved to be eventually
extinct~\cite[Corollary~3]{Bertoin17}, an observation which was
already made by Filippov~\cite{Filippov61} in the context of pure
fragmentations.

Both for homogeneous and for self-similar fragmentations, the additive
martingales~\eqref{eq:A_Gauss} and more precisely their uniform integrability
have turned out to be of greatest importance in the study of asymptotic
behaviors. We stress that sufficient conditions to this uniform
integrability appear less easily for growth-fragmentations, as they non
longer take values in the space of mass-partitions~$\MP$.

Our work is organized in two independent parts. In \cref{sec:A_Riemann}, we
deal with the homogeneous case $\alpha=0$ in the slightly more general
setting of compensated fragmentations~\cite{Bertoin16}. With the help
of a well-known theorem due to Biggins~\cite{Biggins92} and by adapting
arguments of Bertoin and Rouault~\cite{Rouault05}, we prove the uniform
convergence of additive martingales from which, in the realm of
branching random walks, we infer precise estimates for the empirical
measure of the fragments and the asymptotic behavior of the largest one.
This part can be viewed as an application to the study of extremal
statistics in certain branching random walks, see e.g.\ the recent
developments by A\"id\'{e}kon~\cite{Aidekon13}, A\"id\'{e}kon
et al.~\cite{Aidekon12}, Arguin et al.~\cite{Arguin12} and Hu
et al.~\cite{Hu09}. The self-similar case is considered in \cref{sec:A_Poincare}
within the framework of~\cite{Bertoin17}. Relying on recent results
in~\cite{Curien18} and in particular on the uniform integrability of the
Malthusian martingale, we establish for $\alpha>0$ the convergence in
probability of the empirical measure of the fragments and that of the
largest fragment. In a concluding section we also address the
convergence of another empirical measure where fragments are stopped as
soon as they become smaller than a vanishing threshold.

\section{Compensated fragmentations}\label{sec:A_Riemann}
\subsection{Prerequisites}
Recall the space of mass-partitions~$\MP$ defined in~\eqref{eq:A_Newton}
and denote by~$\MP_1$ the subspace of mass-partitions
$\bs p\defeq(p_1,0,\ldots)\in\MP$ having only one single fragment
$p_1\in\oc{0}{1}$. A compensated fragmentation process
$\bs Z(t)\defeq(Z_1(t),Z_2(t),\ldots),\,t\ge0,$ is a stochastic process
whose distribution is characterized by a triple $(\sigma^2,c,\nu)$
where $\sigma^2\ge0$ is a diffusion coefficient, $c\in\RR$ is a growth
rate, and~$\nu$ is a nontrivial measure on
$\MP\setminus\{(1,0,\ldots)\}$ such that \eqref{eq:A_Euler} holds. It can
be seen as the decreasing rearrangement of the exponential of the atoms
of a branching process in continuous time. Namely, the process giving
the empirical measure of the logarithms of the fragments at time~$t$,
\begin{equation*}
  \BRW^t\defeq\sum_{i=1}^\infty\delta_{\log Z_i(t)},
\end{equation*}
is called a \emph{branching L\'{e}vy process} in~\cite{Bertoin16}, to which
we refer for background. In the basic case where
$\nu(\MP\setminus\MP_1)<\infty$, i.e.\ the fragmentation rates are finite,
$\BRW$ is a generalization of the branching random walk in
continuous time introduced by Uchiyama~\cite{Uchiyama82}: more
precisely, $\BRW$ is a branching particle system in which each atom,
during its lifetime, is allowed to move in~$\RR$ independently of the
other atoms and according to the dynamics of a spectrally negative L\'{e}vy
process%
\footnote{That is a c\`{a}dl\`{a}g stochastic process with stationary and
independent increments which has only negative jumps. The results of
this section could be quite straightforwardly adapted to also handle
positive jumps in the particle motions; we shall however not do so as
this would burden the expository and was anyway not considered
in~\cite{Bertoin16}.}~$\eta$ with Laplace transform
\begin{flalign*}
&&&&&&E\!\left[\exp\bigl(q\eta(t)\bigr)\right]&=\,\exp\bigl(t\psi(q)\bigr),&&t\ge0,\ q\ge0,&&
\end{flalign*}
where under~\eqref{eq:A_Euler} the Laplace exponent%
\footnote{Formula~\eqref{eq:A_Lagrange} is designed in such a way that if
  $\sigma^2=0$, $c=0$, and
  $D\defeq\int_{\MP}(1-p_1)\,\nu(\dd\bs p)<\infty$, then~$\bs Z$
  simply is a pure homogeneous fragmentation~$\bs X$ with dislocation
  measure~$\nu$ affected by a dilatation with coefficient~$D$,
  i.e.\ $\bs Z(t)=e^{Dt}\bs X(t),\,t\ge0$. In this case~$\eta$ is a
  compound Poisson process with jump measure
  $(\log p_1)\,\nu_{\mid\MP_1}(\dd\bs p)$ and drift~$D$, but we stress
  that~$\psi(q)<\infty$ holds in greater generality, namely
  under~\eqref{eq:A_Euler} and $\nu(\MP\setminus\MP_1)<\infty$. See
  \cite{Bertoin16} for details.}
\begin{equation}
  \psi(q)\defeq\frac12\sigma^2q^2+\left(c+\int_{\MP\setminus\MP_1}
    (1-p_1)\,\nu(\dd\bs p)\right)\!q
      +\int_{\MP_1}\bigl(p_1^q-1+q(1-p_1)\bigr)\,%
    \nu(\dd\bs p)\label{eq:A_Lagrange}
\end{equation}
is finite for all $q\ge0$.
In words, when $\nu(\MP\setminus\MP_1)<\infty$, the system can be
described as follows. It starts at the origin of space and time with a
single particle which evolves like~$\eta$. Each particle dies after a
random exponential time with intensity $\nu(\MP\setminus\MP_1)$, giving
birth to a random family of children $(\eta_1,\eta_2,\ldots)$ whose
initial position $(\Delta a_1,\Delta a_2,\ldots)$ relative to the mother
particle at its death is such that
$(e^{\Delta a_1},e^{\Delta a_2},\ldots)$ has the conditional
distribution $\nu(\,\cdot\mid\MP\setminus\MP_1)$.

\noindent
In the general situation where the dislocation rate
$\nu(\MP\setminus\MP_1)$ may be infinite, the construction is achieved
by approximation from compensated fragmentations with finite
dislocation rates, using a monotonicity argument (see
\cite[Lemma~3]{Bertoin16} recalled in the proof of \cref{pro:A_Euclid}
below).

Let us denote by
\begin{equation}
  \mu(\dd x)\defeq\EE[\BRW^1(\dd x)]
  \label{eq:A_Hilbert}
\end{equation}
the mean intensity of the point process~$\BRW^1$, so that
\begin{equation*}
  m(q)\defeq\int e^{qx}\mu(\dd x),\qquad q\ge0,
\end{equation*}
is the Laplace transform of~$\mu$. An important fact
(cf.\ \cite[Theorem~1]{Bertoin16}) is that, for every $t\ge0$ and every
$q\ge0$,
\begin{equation}
  m(q)^t\,
    =\,\EE\!\left[\sum_{i=1}^\infty Z_i^q(t)\right]
    =\,\exp\bigl(t\cf(q)\bigr),%
  \label{eq:A_Leibniz}
\end{equation}
where
\begin{equation*}
  \cf(q)\defeq\frac12\sigma^2q^2+cq+\int_\MP\left(\sum_{i=1}^\infty
    p_i^q-1+q(1-p_1)\right)\nu(\dd\bs p)
\end{equation*}
defines a convex function
$\cf\colon\co{0}{\infty}\to\oc{-\infty}{\infty}$. We mention that under
$\nu(\MP\setminus\MP_1)<\infty$, there is the identity
\begin{equation}
  \cf(q)=\psi(q)+\int_{\MP\setminus\MP_1}\!\left(\sum_{i=1}^\infty
    p_i^q-1\right)\nu(\dd\bs p),\qquad q\ge0.%
  \label{eq:A_Grothendieck}
\end{equation}
As we shall explain in the forthcoming \cref{lem:A_Turing}, the first summand
describes the motion of a particle, while the second outlines the
branching mechanism. In better words, $\cf$ is merely the log-Laplace
transform of the cloud of particles at first generation (i.e.\ after the
first branching event), which is a key feature of branching random
walks.

Since under~\eqref{eq:A_Euler},
\begin{equation*}
  p_1^q-1+q(1-p_1)=O\bigl((1-p_1)^2\bigr)
\end{equation*}
is integrable with respect to~$\nu$, we easily observe that, if we set
\begin{equation*}
  \ubar q\defeq\inf\,\bigl\{q\ge0\colon\cf(q)<\infty\bigr\}
    =\inf\left\{q\ge0\colon\int_{\MP\setminus\MP_1}\sum_{i=2}^\infty
      p_i^q\,\nu(\dd\bs p)<\infty\right\}\!,%
\end{equation*}
then~$\cf$ takes finite values and is analytic on the open
interval~$\oo{\ubar q}{\infty}$. Note that~\eqref{eq:A_Euler} also implies
$\cf(2)<\infty$, so $\ubar q\le2$. Let us introduce the subspace
\begin{equation*}
  \DSS[q]\defeq\left\{\bs z\defeq(z_1,z_2,\ldots)\colon z_1\ge z_2
    \ge\cdots\ge0\enspace\text{and}\enspace%
      \sum_{i=1}^{\infty}z_i^q<\infty\right\}
\end{equation*}
of the space~$\ell^q$ of $q$-summable sequences endowed with the
distance $\|\bs z-\bs z'\|_{\ell^q}^q\defeq\sum_{i=1}^\infty|z_i-z'_i|^q$.
We also denote $\DSS[\infty]$ the space of bounded, non-increasing
sequences of nonnegative real numbers endowed with the uniform norm
$\|\cdot\|_{\ell^\infty}$. We see by~\eqref{eq:A_Leibniz} that the
compensated fragmentation $\bs Z\defeq(\bs Z(t),\,t\ge0)$ is a
$\DSS[q]$-valued process for every $q\in\oc{\ubar q}{\infty}$, and in
particular for $q=2$. Further if $\bs z\defeq(z_1,z_2,\ldots)$ is
in~$\DSS$ and $\bs Z^{[1]},\bs Z^{[2]},\ldots$ are independent copies
of~$\bs Z$, then the process of the family
$(z_jZ^{[j]}_i(t),\,i,j\in\NN),\,t\ge0,$ rearranged in the non-increasing
order is again in~$\DSS$, and we denote its distribution by~$\PP_{\bs z}$.
It has been proved in~\cite{Bertoin16} that
$(\bs Z,\,(\PP_{\bs z})_{\bs z\in\DSS})$ is a Markov process which fulfills
the so called \emph{branching property}: for all $s\ge0$, the conditional
law of $(\bs Z(t+s))_{t\ge0}$ given $(\bs Z(r))_{0\le r\le s}$ is
$\PP_{\bs z}$, where $\bs z=\bs Z(s)$. Without loss of generality we shall
assume in the sequel that the fragmentation starts with a single mass
with unit size, i.e.\ $\PP\defeq\PP_{(1,0,\ldots)}$.

Equation~\eqref{eq:A_Leibniz} and the branching property yield an important
family of additive martingales. Namely, the $\RR$-valued process
\begin{equation}
  M(t;q)\defeq\exp\bigl(-t\cf(q)\bigr)\sum_{i=1}^\infty Z_i^q(t),\qquad
    t\ge0,\label{eq:A_Fermat}
\end{equation}
is a martingale for every $q\in\oo{\ubar q}{\infty}$. As a first
consequence \cite[Proposition~2]{Bertoin16}, the compensated
fragmentation~$\bs Z$ possesses a c\`{a}dl\`{a}g version in~$\DSS$, that is a
version in the Skorokhod space $D(\co{0}{\infty},\DSS)$ of right
continuous with left limits, $\DSS$-valued functions. Working with such
a version from now on, $\bs Z$ has actually c\`{a}dl\`{a}g paths in~$\DSS[q]$
for every $q\in\oc{\ubar q}{\infty}$.

\begin{proposition}\label{pro:A_Galois}
  Almost surely, for every $q\in\oc{\ubar q}{\infty}$, $\bs Z$ has c\`{a}dl\`{a}g
  paths in~$\DSS[q]$.
\end{proposition}
\begin{proof}
  Recall that $\|\cdot\|_{\ell^{q'}}\le\|\cdot\|_{\ell^q}$ whenever
  $q\le q'\le\infty$. Since~$\bs Z$ has c\`{a}dl\`{a}g paths in~$\DSS$, it has
  in particular c\`{a}dl\`{a}g paths in~$\DSS[\infty]$. Let
  $(q_k,\,k\in\NN)$ be a sequence decreasing to~$\ubar q$, and define
\begin{equation*}
T^{(k)}_m
    \defeq\inf\,\bigl\{t\ge0\colon M(t;q_k)>m\bigr\}
    =\inf\,\Bigl\{t\ge0\colon \|\bs Z(t)\|_{\ell^{q_k}}^{q_k}
      >m\,e^{t\cf(q_k)}\Bigr\}
\end{equation*}
  for $k,m\in\NN$. Applying Doob's maximal inequality to the
  martingale~\eqref{eq:A_Fermat} we have that almost surely,
  for all $k\in\NN$, $T^{(k)}_m\uparrow\infty$ as $m\to\infty$.
  Thus, almost surely, for every $q\in\oc{\ubar q}{\infty}$ and
  $T\ge0$ we can find a $k\in\NN$ such that $\ubar q<q_k<q$ and then a
  $m\in\NN$ such that $T<T^{(k)}_m$, whence
\begin{align*}
  \|\bs Z(t)-\bs Z(s)\|_{\ell^q}^q
  &\le\|\bs Z(t)-\bs Z(s)\|_{\ell^\infty}^{q-q_k}\,
    \|\bs Z(t)-\bs Z(s)\|_{\ell^{q_k}}^{q_k}\\[.4em]
  &\le m\,2^{1+q_k}\!\left(1+e^{T_m^{(k)}\cf(q_k)}\right)
    \|\bs Z(t)-\bs Z(s)\|_{\ell^\infty}^{q-q_k}
\end{align*}
  for all $0\le s,t<T$. The fact that~$\bs Z$ has c\`{a}dl\`{a}g paths
  in~$\DSS[\infty]$ completes the proof.
\end{proof}

We first would like to extend to the compensated fragmentation~$\bs Z$
the asymptotic results obtained by Bertoin and
Rouault~\cite{Rouault04,Rouault05} for pure homogeneous fragmentations.
They strongly rely on the work of \cite{Biggins77,Biggins92} about the
uniform integrability of additive martingales. Essentially, the
martingales $(M(t;q))_{t\ge0}$ will be uniformly integrable if
$q\cf'(q)-\cf(q)<0$ and $M(1;q)\in\Leb^\gamma(\PP)$ for some
$\gamma>1$. With this in mind, let us introduce
\begin{equation*}\bar q\defeq\sup\,\Bigl\{q>\ubar q\colon
  q\cf'(q)-\cf(q)<0\Bigr\}.
\end{equation*}
First note that $\bar q<\infty$, because
\begin{equation*}
q\cf'(q)-\cf(q)=\frac12\sigma^2q^2+\int_{\MP}
  \bigl(1-p_1^q(1-\log p_1^q)\bigr)\,\nu(\dd\bs p)
  -\int_{\MP\setminus\MP_1}
    \sum_{i=2}^{\infty}p_i^q\bigl(1-\log p_i^q\bigr)\,\nu(\dd\bs p),
\end{equation*}
which, by Fatou's lemma, is at least $\nu(\MP)$ as $q\to\infty$.
Second, we have $\bar q>\ubar q$ as soon as $q\cf'(q)-\cf(q)<0$ for some
$q\ge0$ such that $\cf(q)<\infty$ (e.g.\ for $q=2$), which is
realized when
\begin{equation*}
\int_{\MP\setminus\MP_1}\sum_{i=2}^{\infty}p_i^q\bigl(1-\log p_i^q
    \bigr)\,\nu(\dd\bs p)
  \ge\left(\frac12\sigma^2+\int_{\MP}(1-p_1)^2\,\nu(\dd\bs p)
    \right)\!q^2.
\end{equation*}
We distinguish two different regimes for the function $q\mapsto\cf(q)/q$:
\begin{lemma}\label{lem:A_Neumann}
The function $q\mapsto\cf(q)/q$ is decreasing on~$\oo{\ubar q}{\bar q}$
and increasing on~$\oo{\bar q}{\infty}$.
Further, $\bar q\cf'(\bar q)=\cf(\bar q)$ when $\bar q>\ubar q$.
\end{lemma}
In the context of branching random walks, the value $\cf'(\bar q)$ is the
asymptotic velocity of the maximal displacement~$\log Z_1(t)$; see
\cref{fig:A_Abel} and \cref{sec:A_Weierstrass}.

\begin{figure}[htbp]
  \begin{subfigure}[b]{0.34\textwidth}
    \centering
  \begin{tikzpicture}
    \draw[->] (-.25,0) -- (3.8,0);
    \draw[->] (0,-.25) -- (0,3.25) node [left] {$\cf(q)$};
    \draw plot [domain=4/13:3/2] (\x, {1/\x+(2/3)^2*\x-1});
    \draw plot [domain=3/2:13/4] (\x, {(\x-3/2)^(3/2)+1/3});
    \draw[dotted] plot [domain=0:3.75]
      (\x, {0.2211539327017409*\x});
    \draw[->] (1.033533639425775,0.2285700289385531)
      -- (2.009941082307065,0.4445063748510008)
        node [below right=-3pt] {\footnotesize $\cf'(\bar q)$};
    \draw (1.52173736086642,-.05) -- (1.52173736086642,0.05);
    \draw[dotted]
      (1.52173736086642,0) node [below] {$\bar q\vphantom{\ubar q}$}
         -- (1.52173736086642,0.336538201894777);
    \draw[dotted] (0.9*4/13,0) node [below]
      {$\ubar q\vphantom{\bar q}$} -- (0.9*4/13,2.75);
  \end{tikzpicture}

  {\bfseries\sffamily a.} Positive velocity.
  \end{subfigure}%
  \begin{subfigure}[b]{0.32\textwidth}
    \centering
  \begin{tikzpicture}
    \draw[->] (-.25,0) -- (3.3,0);
    \draw[->] (0,-.25) -- (0,3.25) node [left] {$\cf(q)$};
    \draw plot[domain=(4/13:3/2)] (\x, {(\x-3/2)^2)});
    \draw plot[domain=3/2:3] (\x, {(\x-3/2)^3});
    \draw[dotted] (4/13,0) node [below] {$\ubar q\vphantom{\bar q}$}
      -- (4/13,961/676) node {$\scriptscriptstyle\bullet$};
    \draw (3/2,0) node [below] {$\bar q\vphantom{\ubar q}$};
    \draw (4/13,-0.05) -- (4/13,0.05);
    \draw (3/2,-0.05) -- (3/2,0.05);
  \end{tikzpicture}

  {\bfseries\sffamily b.} Positive velocity.
  \end{subfigure}%
  \begin{subfigure}[b]{0.33\textwidth}
    \centering
  \begin{tikzpicture}
    \draw[->] (-.25,0) -- (3.5,0);
    \draw[->] (0,-2.25) -- (0,1.25) node [left] {$\cf(q)$};
    \draw plot [domain=(4/13:21/13)] (\x,{(\x-21/13)^2-2});
    \draw plot [domain=(21/13:7/2)] (\x,{(\x-21/13)^(7/4)-2});
    \draw[dotted] (0,0) -- (1.25,-2.08675162209765);
    \draw[dotted] (4/13,-0.2899408284023668)
      node {$\scriptscriptstyle\bullet$}
       -- (4/13,0) node [above] {$\ubar q\vphantom{\bar q}$};
    \draw[dotted] (0.7806839665455554,0)
      node [above] {$\bar q\vphantom{\ubar q}$} --
      (0.7806839665455554,-1.303274826827652);
    \draw (4/13,-0.05) -- (4/13,0.05);
    \draw (0.7806839665455554,-0.05) -- (0.7806839665455554,0.05);
    \draw[->] (0.5237461399564401,-0.874342485697187)
      -- (1.037621793134671,-1.732207167958117) node[left]
      {\footnotesize$\cf'(\bar q)$};
   \draw (0,-2) node [below] {$\vphantom{\bar{\ubar q}}$};
  \end{tikzpicture}

  {\bfseries\sffamily c.} Positive velocity.
  \end{subfigure}
  \caption{The cumulant function~$\cf$,
    the points~$\ubar q,\bar q$, and the velocity
    $\cf'(\bar q)=\cf(\bar q)/\bar q$.}%
  \label{fig:A_Abel}
\end{figure}

\begin{proof}
On the one hand,
\begin{equation*}
  \frac{\dd}{\dd q}\left[\frac{\cf(q)}q\right]
    =\,\frac{q\cf'(q)-\cf(q)}{q^2}.
\end{equation*}
On the other hand, the map $q\mapsto q\cf'(q)-\cf(q)$ is
increasing on~$\oo{\ubar q}{\infty}$ since~$\cf$ is convex, so it
has at most one sign change, occurring at~$\bar q$ if $\bar q>\ubar q$.
\end{proof}

Our main result provides sufficient conditions for the convergence of
the martingales $(M(t;q))_{t\ge0}$ uniformly in
$q\in\oo{\ubar q}{\bar q}$, both almost surely and in $\Leb^1(\PP)$.
Most of the coming section is devoted to a precise statement and a
proof. As consequences, we ascertain the convergence of a rescaled
version of the empirical measure~$\BRW^t$ and in \cref{sec:A_Weierstrass} we expand
on the asymptotic behavior of the largest fragment. One last application
is exposed in \cref{sec:A_Descartes}.

\subsection{Uniform convergence of the additive martingales}
In the remaining of \cref{sec:A_Riemann} we will make, in addition
to~\eqref{eq:A_Euler}, the assumption that the dislocation measure~$\nu$
fulfills
\begin{gather}
  \cf(0)\in\oc{0}{\infty},\label{eq:A_Dirichlet}
\intertext{and, for all $\ubar q<q<1$,}
  \nu_{|\MP\setminus\MP_1}\!\left(\sum_{i=1}^\infty p_i^q<1\right)
    <\,\infty.\label{eq:A_Ramanujan}
\end{gather}
Condition~\eqref{eq:A_Dirichlet} holds e.g.\ when $\nu(p_2>0)>\nu(p_1=0)$
and merely rephrases that the mean number $\mu(\RR)=m(0)$ of offspring of
particles is greater than~$1$, i.e.\ the branching process~$\BRW$ is
supercritical. This implies that the \emph{non-extinction} event
$\{\forall t\ge0,\ Z_1(t)>0\}$ occurs with positive probability.
Condition~\eqref{eq:A_Ramanujan} is just a minor technical requirement for the
possible values $q<1$ and is fulfilled in many situations: when
$\ubar q\ge1$, when $\nu(\MP\setminus\MP_1)<\infty$, or more
importantly when the measure~$\nu$ is \emph{conservative}, i.e.\
$\sum_{i\ge1}p_i=1$ for $\nu$-almost every $\bs p\in\MP$. Observe also
that in the conservative case, $\ubar q<1$ is possible only
if~\eqref{eq:A_Archimedes} holds, i.e.~$\bs Z$ is essentially a dilated pure
fragmentation.

\medskip
We may now state:
\begin{theorem}\label{thm:A_Jacobi}
  Suppose~\eqref{eq:A_Dirichlet} and~\eqref{eq:A_Ramanujan}. Then the following
  assertions hold almost surely:
\begin{enumerate}[label=(\roman*)]
  \item On~$\oo{\ubar q}{\bar q}$, $M(t;\cdot)$ converges locally
    uniformly as $t\to\infty$. More precisely, there exists a random
    continuous function
    $M(\infty;\cdot)\colon\oo{\ubar q}{\bar q}\to\co{0}{\infty}$ such
    that, for any compact subset~$K$ of~$\oo{\ubar q}{\bar q}$,
  \begin{equation*}
    \lim_{t\to\infty}\sup_{q\in K}|M(t;q)-M(\infty;q)|=0,
  \end{equation*}
    and this convergence also holds in mean.
    Furthermore for every $q\in\oo{\ubar q}{\bar q}$, $M(\infty;q)>0$
    conditionally on non-extinction.\label{thm:A_Jacobi1}
  \item For every $q\in\co{\bar q}{\infty}$,\label{thm:A_Jacobi2}
    \begin{equation*}
      \lim_{t\to\infty}M(t;q)=0.
    \end{equation*}
\end{enumerate}
\end{theorem}
As a first important consequence, we derive uniform estimates for the
empirical measure of the fragments, which echo those determined by
Bertoin and Rouault~\cite[Corollary~3]{Rouault05}. We shall assume here
that the mean intensity measure~$\mu$ in~\eqref{eq:A_Hilbert} is
\emph{non-lattice}, in that it is not supported on~$r\ZZ+s$ for any
$r>0,\,s\in\RR$.
\begin{corollary}\label{cor:A_Cantor}
  Suppose~\eqref{eq:A_Dirichlet}, \eqref{eq:A_Ramanujan}, and $\mu$ non-lattice.
  Then for any Riemann integrable function
  $f\colon\oo{0}{\infty}\to\RR$
  with compact support and for all compact subset~$K$
  of~$\oo{\ubar q}{\bar q}$,
\begin{equation*}
  \lim_{t\to\infty}\sqrt t\,e^{-\bigl(\cf(q)-q\cf'(q)\bigr)t}
    \sum_{i=1}^\infty f\!\left(Z_i(t)\,e^{-\cf'(q)t}\right)
  =\,\frac{M(\infty;q)}{\sqrt{2\pi\cf''(q)}}\int_0^\infty
    \frac{f(y)}{y^{q+1}}\,\dd y
\end{equation*}
  uniformly in $q\in K$, almost surely.
\end{corollary}
\begin{remark}
We stress that condition~\eqref{eq:A_Ramanujan} is unnecessary if we only deal
with $q\ge1$. In particular it may be removed from the above statements
provided that we replace~$\ubar q$ by $\ubar q\vee 1$.
\end{remark}

Before proving these two results, let us give a quick summary on the
sizes of particles in a compensated fragmentation. On the one hand,
it is easy (see e.g.\ \cite[Corollary~1.4]{Bertoin06}) to derive from
\hyperref[thm:A_Jacobi1]{\cref*{thm:A_Jacobi}.(i)}
that in the first order, the largest particle $Z_1(t)$ evolves like
$e^{\cf'(\bar q)t}$ as $t\to\infty$, and we will have a look at the
second and third asymptotic orders in~\cref{sec:A_Weierstrass}. On the other
hand, \cref{cor:A_Cantor} provides the local density of particles at
intermediate scales: if $\cf'(\ubar q)<a<\cf'(\bar q)$ and
$\cf^*(a)\defeq\cf(q)-q\cf'(q)$ for $\cf'(q)\defeq a$, then
\begin{equation*}
  \lim_{t\to\infty}\frac1t\log
  \#\!\left\{i\in\NN\colon e^{at-\eps}\le Z_i(t)\le
    e^{at+\eps}\right\}=\,\cf^*(a)
\end{equation*}
for every $\eps>0$, almost surely (just take
$f(x)\defeq\II_{\cc{-\eps}{\eps}}(\log x)$ above).
Lastly, we shall observe in \cref{sec:A_Descartes} that fragments at
untypical levels $a>\cf'(\bar q)$ appear with a probability that is
roughly of the same order as their expected number
(\hyperref[cor:A_Weyl2]{\cref*{cor:A_Weyl}.(ii)}).

\cref{thm:A_Jacobi} is essentially a version of a theorem of
Biggins~\cite{Biggins92} in the context of compensated fragmentations.
In this respect, one important requirement to derive
part~\labelcref{thm:A_Jacobi1} is that $\EE[M(1;q)^\gamma]<\infty$ for some
$\gamma>1$. We start with a lemma controlling the finiteness of
\begin{equation*}
W^\gamma_{\nu,q}\defeq\int_{\MP\setminus\MP_1}
\left|1-\sum_{i=1}^\infty p_i^q\right|^\gamma\nu(\dd\bs p).
\end{equation*}

\begin{lemma}\label{lem:A_Cayley}
  Let $q>\ubar q$ and suppose either~\eqref{eq:A_Ramanujan} or $q\ge1$. Then
  $W^\gamma_{\nu,q}<\infty$ for some $\gamma\in\oc{1}{2}$.
\end{lemma}
\begin{proof}
  Suppose first $q\ge1$. Then for $\gamma\defeq2$ and for all
  $\bs p\in\MP$,
\begin{equation*}
  0\le\left(1-\sum_{i=1}^\infty p_i^q\right)^{\!2}\le(1-p_1^q)^2
    \le q^2\,(1-p_1)^2
\end{equation*}
  (the last inequality resulting from the convexity of
   $x\mapsto x^q$), so $W^2_{\nu,q}<\infty$ by~\eqref{eq:A_Euler}.
  Suppose now $q<1$. Then
\begin{align*}
  W^\gamma_{\nu,q}
    \le\nu_{|\MP\setminus\MP_1}\!\left(\sum_{i=1}^\infty p_i^q<1\right)
      +\int_{\MP\setminus\MP_1}\II_{\left\{\sum\limits_{i=1}^\infty p_i^q
        \ge1\right\}}\!\left(\sum_{i=1}^\infty p_i^q-1\right)^{\!\gamma}
        \nu(\dd\bs p).
\end{align*}
  Under~\eqref{eq:A_Ramanujan}, $W^\gamma_{\nu,q}$ is finite as soon as the
  latter integral is finite. But by Jensen's inequality, the integrand
  is bounded from above by
\begin{equation*}
  \left(\sum_{i=2}^\infty p_ip_i^{q-1}\right)^{\!\gamma}
    \le\,\sum_{i=2}^\infty p_i^{1+\gamma(q-1)},
\end{equation*}
  which is $\nu$-integrable if
  $1+\gamma(q-1)=q-(\gamma-1)(1-q)>\ubar q$, i.e.\ provided
  $\gamma\in\oc{1}{2}$ is close enough to~$1$.
\end{proof}

We then derive an upper bound for $\EE[M(t;q)^\gamma]$ in terms of
$W^\gamma_{\nu,q}$:
\begin{lemma}\label{lem:A_Turing}
  Suppose that $\nu(\MP\setminus\MP_1)<\infty$. Then for every
  $q\in\oo{\ubar q}{\infty}$, $\gamma\in\oc{1}{2}$ and $t\ge0$,
\begin{equation}
  \EE[M(t;q)^\gamma]
    \le\,\mathfrak{c}_\gamma\,W^\gamma_{\nu,q}\,f\bigl(t,\psi(\gamma q)
      -\gamma\psi(q),\cf(\gamma q)-\gamma\cf(q)\bigr),%
  \label{eq:A_Noether}
\end{equation}
  where $\psi$ is given by~\eqref{eq:A_Lagrange},
  $f(t,x,y)\defeq(e^{tx}-e^{ty})/(x-y)$, and $\mathfrak{c}_\gamma$ is
  a finite constant depending only on~$\gamma$.
\end{lemma}
\begin{proof}
  Lemma~2 in~\cite{Bertoin16} states that the branching L\'{e}vy
  process~$\mathcal{Z}$ can be obtained by superposing independent
  spatial L\'{e}vy motions to a ``steady'' branching random walk.
  Specifically, for each $t\ge0$,
\begin{equation*}
\bs Z(t)
    \eqlaw\left(e^{\beta_1}X_1(t),e^{\beta_2}X_2(t),\ldots\right)\!,
\end{equation*}
  where $\bs X(t)\defeq(X_1(t),X_2(t),\ldots)$ are the atoms at time $t$
  of a homogeneous fragmentation~$\bs X$ with dislocation measure
  $\nu_{|\MP\setminus\MP_1}$ and $(\beta_i)_{i\in\NN}$ is an independent
  sequence of random variables with Laplace transform
  $\EE[\exp(q\beta_i)]=\exp(t\psi(q)),\,q\ge0$.
  Applying Jensen's inequality and conditioning on~$\bs X(t)$ produce
\begin{align}
  \EE\!\left[\left(\sum_{i=1}^\infty Z_i^q(t)\right)^{\!\gamma\,}\right]
  &=\,\EE\!\left[\left(\sum_{j=1}^\infty X_j^q(t)\right)^{\!\gamma}
    \left(\sum\limits_{i=1}^\infty e^{q\beta_i}
      \frac{X_i^q(t)}{\sum_jX_j^q(t)}\right)^{\!\gamma\,}\right]\notag%
      \\[.4em]
  &\le\,\EE\!\left[\left(\sum_{j=1}^\infty X_j^q(t)\right)^{\!\gamma-1}
      \sum_{i=1}^\infty e^{\gamma q\beta_i}
        X_i^q(t)\right]\notag\\[.4em]
  &=\,\exp\bigl(t\psi(\gamma q)\bigr)\,\EE\!\left[\left(\sum_{i=1}^\infty X^q_i(t)
      \right)^{\!\gamma\,}\right]\!.\label{pro:A_Turing1}
\end{align}
  We now recall from~\cite{Bertoin03} (see the proof of its Theorem~2)
  how to estimate the latter expectation. Denoting
\begin{equation*}
\phi(q)
    \defeq\int_{\MP\setminus\MP_1}\!\left(\sum_{i=1}^\infty p_i^q-1
      \right)\nu(\dd\bs p)<\infty,
\end{equation*}
  the process
\begin{equation}
  N(t;q)
    \defeq\exp\bigl(-t\phi(q)\bigr)
      \sum_{i=1}^\infty X_i^q(t),\qquad t\ge0,%
  \label{eq:A_Pythagoras}
\end{equation}
  is a purely discontinuous martingale. It is then deduced from an
  inequality of Burkholder--Davis--Gundy that
  \begin{equation*}
    \EE[N(t;q)^\gamma]\le\mathfrak{c}_\gamma\,\EE[V^\gamma(t;q)],
  \end{equation*}
  where $\mathfrak{c}_\gamma<\infty$ is some constant, and~$V^\gamma$
  is the $\gamma$-variation process of~$N$:
  \begin{equation*}
    V^\gamma(t;q)\defeq\sum_{0<s\le t}|N(s;q)-N(s-;q)|^\gamma.
  \end{equation*}
  Since in this setting
  \begin{equation*}
    |N(s;q)-N(s-;q)|^\gamma=\exp\bigl(-s\gamma\phi(q)\bigr)
      X_k^{\gamma q}(s-)\left|1-\sum_{i=1}^\infty p_i^q\right|^\gamma,
  \end{equation*}
  $(s,\bs p,k)\in\oc{0}{t}\times\MP\times\NN$ being the atoms of a Poisson
  random measure
  with intensity $\dd t\otimes\nu_{|\MP\setminus\MP_1}\otimes\sharp$,
  it follows that~$V^\gamma$ has predictable compensator
\begin{equation*}
  \left(\int_{\MP\setminus\MP_1}\left|\sum_{i=1}^\infty
     1-p_i^q\right|^\gamma\nu(\dd\bs p)\right)
   \int_0^t\exp\bigl(-s\gamma\phi(q)\bigr)\sum_{i=1}^\infty
     X_i^{\gamma q}(s)\,\dd s,
\end{equation*}
and therefore
\begin{equation}
  \EE[N(t;q)^\gamma]
    \le \mathfrak{c}_\gamma\,W^\gamma_{\nu,q}\,f\bigl(t,0,\phi(\gamma q)
      -\gamma\phi(q)\bigr).%
  \label{pro:A_Turing2}
\end{equation}
  Now recall~\eqref{eq:A_Fermat}, \eqref{eq:A_Pythagoras} and the identity
  $\phi(q)+\psi(q)=\cf(q)$ already observed in~\eqref{eq:A_Grothendieck}.
  Multiplying~\eqref{pro:A_Turing1} by
  $e^{-t\gamma\cf(q)}$ and then reporting the bound~\eqref{pro:A_Turing2},
  we end up with
\begin{equation*}
  \EE[M(t;q)^\gamma]
    \le \mathfrak{c}_\gamma\,W^\gamma_{\nu,q}\,f\bigl(t,\psi(\gamma q)
      -\gamma\psi(q),\cf(\gamma q)-\gamma\cf(q)\bigr),
\end{equation*}
  as desired.
\end{proof}

Putting the previous results together now yields:
\begin{proposition}\label{pro:A_Euclid}
  Let $q>\ubar q$ and suppose either~\eqref{eq:A_Ramanujan} or $q\ge1$.
  Then there exists $\gamma\in\oc{1}{2}$ such that
  $M(t;q)\in\Leb^\gamma(\PP)$ for all $t\ge0$.
\end{proposition}
\begin{proof}
  By \cref{lem:A_Cayley} we can choose $\gamma\in\oc{1}{2}$ such that
  $W^\gamma_{\nu,q}<\infty$. Let us first assume
  $\nu(\MP\setminus\MP_1)<\infty$, so that we may
  apply~\cref{lem:A_Turing}. Note that
\begin{equation*}
\psi(\gamma q)-\gamma\psi(q)=\frac12\sigma^2(\gamma-1)\gamma q^2
    +\int_{\MP_1}(p_1^{\gamma q}-\gamma p_1^q+\gamma-1)\,\nu(\dd\bs p),
\end{equation*}
  with
\begin{equation*}
  0\le p_1^{\gamma q}-\gamma p_1^q+\gamma-1=O\bigl((1-p_1)^2\bigr).
\end{equation*}
  Then the inequality~\eqref{eq:A_Noether} is
\begin{equation*}
  \EE[M(t;q)^\gamma]
  \le\,\mathfrak{c}_\gamma\,W^\gamma_{\nu,q}\,f\!\left(t,\frac12\sigma^2
    (\gamma-1)\gamma q^2+\int_{\MP_1}(p_1^{\gamma q}-\gamma p_1^q
    +\gamma-1)\,\nu(\dd\bs p),\cf(\gamma q)-\gamma\cf(q)\right)\!,
\end{equation*}
  where~$f$ is a continuous function. This bound is finite for each
  $t\ge0$, and we shall show by approximation that this also holds when
  $\nu(\MP\setminus\MP_1)=\infty$. The measures~$\nu^{(b)}$, images
  of~$\nu$ by the maps
\begin{equation*}
  \bs p\longmapsto\left(p_1,p_2\II_{\{p_2>e^{-b}\}},
    p_3\II_{\{p_3>e^{-b}\}},\ldots\right)\!,\qquad b>0,
\end{equation*}
  define a consistent family of dislocation measures such
  that $\nu^{(b)}(\MP\setminus\MP_1)<\infty$. Thanks to
  \cite[Lemma~3]{Bertoin16} we can consider that~$\bs Z$ arises from the
  inductive limit $\bs Z\defeq\mathop{\lim\!\!\uparrow}\bs Z^{(b)}$ as $b\uparrow\infty$,
  where the $\bs Z^{(b)},\,b>0,$ are suitably embedded compensated
  fragmentations with characteristics $(\sigma^2,c,\nu^{(b)})$. With
  obvious notations, we deduce from the monotone convergence theorem
  that
  $\EE[\|\bs Z^{(b)}(t)\|_{\ell^q}^{\gamma q}]
    \to\EE[\|\bs Z(t)\|_{\ell^q}^{\gamma q}]$
  as $b\to\infty$, and from the dominated convergence theorem that
  $\cf^{(b)}(q)\to\cf(q)$ and
  $W^\gamma_{\nu^{(b)},q}\to W^\gamma_{\nu,q}$
  (working like in the proof of~\cref{lem:A_Cayley}).
  The proof is then completed by Fatou's lemma.
\end{proof}

We are finally ready to tackle the proof of the convergence of the
martingale $(M(t;q))_{t\ge0}$.
\begin{proof}[Proof of~\cref{thm:A_Jacobi}]
  Since $(M(t;q))_{t\ge0}$ is a nonnegative c\`{a}dl\`{a}g martingale, its limit
  $M(\infty;q)$ as $t\to\infty$ exists almost surely.
  If $q\in\co{\bar q}{\infty}$, then $q\cf'(q)-\cf(q)\ge0$ so that
  condition~(3.3) in~\cite{Biggins77} fails, and therefore
  $M(\infty;q)=0$. This proves~\labelcref{thm:A_Jacobi2}.
  For~\labelcref{thm:A_Jacobi1}, we follow the lines
  of~\cite{Rouault05}. From \cref{pro:A_Galois} we know that almost
  surely, for every $t\ge0$ and $(t_n,\,n\in\NN)$ such that
  $t_n\downarrow t$ as $n\to\infty$, the sequence of random
  functions on~$K$
\begin{flalign*}
&&&&&&&&q&\mapsto\bigl(1+Z_1(t_n)\bigr)^{-q}\sum_{i=1}^\infty Z_i^q(t_n),&n\in\NN,%
&&&&
\intertext{which all are non-increasing because of the leading factors
$(1+Z_1(t_n))^{-q}$, converges pointwise to the random continuous
function}
&&&&&&&&q&\mapsto\bigl(1+Z_1(t)\bigr)^{-q}\sum_{i=1}^\infty Z_i^q(t).
\end{flalign*}
  By a classical counterpart of Dini's theorem
  (see e.g.\ \cite[Problem~II.3.127]{Polya72}), the convergence is
  actually uniform in $q\in K$. Multiplying by the continuous function
  $q\mapsto(1+Z_1(t))^q\exp(-t\cf(q))$ and dealing similarly with the
  left limits of~$\bs Z$, we can therefore view $(M(t;\cdot))_{t\ge0}$ as
  a martingale with c\`{a}dl\`{a}g paths in the Banach space $\Co(K,\RR)$ of
  continuous functions on~$K$.

  We now observe that the process
  \begin{equation*}
    \BRW^n=\sum_{i=1}^\infty\delta_{\log Z_i(n)},%
    \qquad n\in\NN,
  \end{equation*}
  is a branching random walk (in discrete time), and check the two
  conditions to apply the results of Biggins
  \cite[Theorems~1 \& 2]{Biggins92}: first, if
  $q\in\oo{\ubar q}{\bar q}$ then by \cref{pro:A_Euclid} we have
  $\EE[M(1;q)^\gamma]<\infty$ for some $\gamma\in\oc{1}{2}$; second,
  using \cref{lem:A_Neumann} we can find $\alpha\in\oc{1}{\gamma}$ such
  that $\alpha q\in\oo{q}{\bar q}$, hence
\begin{equation}
  \frac{m(\alpha q)}{m(q)^{\alpha}}
    =\exp\left\{\alpha q\left(\frac{\cf(\alpha q)}{\alpha q}
      -\frac{\cf(q)}q\right)\right\}<\,1.%
   \label{eq:A_Fibonacci}
\end{equation}
  It thus follows that the $\Co(K,\RR)$-valued discrete-time martingale
  \begin{equation*}
    M(n;\cdot)\colon q\mapsto\exp\bigl(-n\cf(q)\bigr)
      \int e^{qx}\,\BRW^{n}(\dd x),\qquad n\in\NN,
  \end{equation*}
  converges as $n\to\infty$ to a random function
  $M(\infty;\cdot)\in\Co(K,\RR)$, almost surely and in mean. Now, the
  uniform norm $\|\cdot\|$ of $\Co(K,\RR)$ is a convex map and thus for
  any integer $n\ge0$ the process $(\|M(t;\cdot)-M(n;\cdot)\|)_{t\ge n}$
  is a nonnegative submartingale with c\`{a}dl\`{a}g paths. If $t\ge0$ and~$n$
  is chosen so that $n\le t<n+1$, we have in particular
\begin{equation*}
  \EE[\|M(t;\cdot)-M(n;\cdot)\|]\le\EE[\|M(n+1;\cdot)-M(n;\cdot)\|],
\end{equation*}
  and consequently
\begin{equation*}
  \EE[\|M(t;\cdot)-M(\infty;\cdot)\|]
    \le\EE[\|M(n+1;\cdot)-M(n;\cdot)\|]
      +\EE[\|M(n;\cdot)-M(\infty;\cdot)\|].
\end{equation*}
  The convergence in $\Leb^1(\PP)$ of the continuous-time martingale
  $(M(t;\cdot))_{t\ge0}$ then follows from the one in discrete time.
  The almost sure convergence is established by applying Doob's maximal
  inequality and the Borel--Cantelli lemma, like in the proof
  of~\cite{Rouault05}: indeed for every $\eps>0$,
\begin{equation*}
  \PP(\exists t\ge n\colon\|M(t;\cdot)-M(n;\cdot)\|>\eps)\,
  \le\,\eps^{-1}\,\EE[\|M(\infty;\cdot)-M(n;\cdot)\|]\,%
    \xrightarrow[n\to\infty]{}\,0.
\end{equation*}

  We finally deal with the almost sure positivity of the terminal value
  $M(\infty;q)$ conditionally on non-extinction. We derive from the
  branching property at time $n$ that, for every
  $q\in\oo{\ubar q}{\bar q}$,
\begin{equation*}
  \PP\!\left(M(\infty;q)=0\;\big|\;\BRW^n\right)
    =\prod_{z\in\BRW^n}\!\PP_z(M(\infty;q)=0),
\end{equation*}
  where by scaling, the probability $\PP_z(M(\infty;q)=0)$ does actually
  not depend on the initial size~$z$. Hence
  $\rho\defeq\PP(M(\infty;q)=0)=\EE[\rho^{\#(n)}]$,
  where $\#(n)\defeq\BRW^n(\RR)$, the number of particles at time
  $n\in\NN$, defines a supercritical Galton--Watson process. Since
  $\rho<1$ (because $\EE[M(\infty;q)]=\EE[M(0;q)]=1$), $\rho$ is its
  probability of extinction. The two events
  $\{\text{extinction}\}\subseteq\{M(\infty;q)=0\}$ having thus the same
  probability we conclude that they coincide up to a negligible event.
\end{proof}
\begin{remark}
  Under the conditions of~\cref{thm:A_Jacobi} we have also from
  \cite[Theorem~5]{Biggins92} that for each $q\in\oo{\ubar q}{\bar q}$
  and $\alpha\in\oc{1}{\gamma}$ as in~\eqref{eq:A_Fibonacci}, the martingale
  $(M(t;q))_{t\ge0}$ converges in $\Leb^\alpha(\PP)$.
\end{remark}
We close this section with the proof of \cref{cor:A_Cantor}.
\begin{proof}[Proof of \cref{cor:A_Cantor}]
  Let us define the tilted measures
\begin{equation*}
  \BRW_q^t(\dd x)\defeq\frac{e^{qx}}{m(q)^t}\,\BRW^t(\dd x),\quad t\ge0,%
\qquad\text{and}\qquad
  \mu_q(\dd x)
    \defeq\frac{e^{qx}}{m(q)}\,\mu(\dd x).
\end{equation*}
  Using~\eqref{eq:A_Leibniz}, $\mu_q$ is a probability measure with mean
\begin{align*}
  c_q&\defeq m(q)^{-1}\,\EE\!\left[\sum_{i=1}^\infty Z_i^q(1)\log Z_i(1)
    \right]\!=\cf'(q),
\shortintertext{and variance}
  \sigma_q^2&\defeq m(q)^{-1}\,\EE\!\left[\sum_{i=1}^\infty Z_i^q(1)
    \log^2 Z_i(1)\right]\!-c_q^2=\cf''(q).
\end{align*}
  On the one hand, we observe that for every $n\in\NN$,
\begin{equation*}
  e^{-\bigl(\cf(q)-q\cf'(q)\bigr)n}
    \sum_{i=1}^\infty f\!\left(Z_i(n)\,e^{-\cf'(q)n}\right)
  =\int_\RR f(e^x)\,e^{-qx}\,\BRW^n_q(nc_q+\dd x).
\end{equation*}
  On the other hand, by a local limit theorem due to Stone
  \cite[Theorem~2]{Stone67},
\begin{align*}
  &&\sqrt n\,\mu_q^{(\star n)}(nc_q+\dd x)
    \,\approx\,p_q\!\left(\frac{x}{\sqrt n}\right)\dd x,%
    &&&n\to\infty,&&
  \intertext{uniformly for $x\in\RR$ and~$q$ in compact subsets
  of~$\oo{\ubar q}{\infty}$, where $\mu_q^{(\star n)},\,n\in\NN,$ is
  the $n$\textsuperscript{th} convolution of~$\mu_q$ with itself and
  $p_q(x)$ denotes the Gaussian density with mean~$0$ and
  variance~$\sigma^2_q$. Thanks to the uniform convergence in
  \cref{thm:A_Jacobi}, this in terms of the branching random walk
  translates into}
  &&\sqrt n\,\BRW_q^n(nc_q+\dd x)
    \,\approx\,M(\infty;q)\,p_q\!\left(\frac{x}{\sqrt n}\right)\dd x,%
    &&&n\to\infty,
\end{align*}
  uniformly for $x\in\RR$ and~$q$ in compact subsets
  of~$\oo{\ubar q}{\bar q}$, almost surely. The corollary then results
  from a Riemann sum argument. We leave details and refer the
  interested reader to Corollary~4 in~\cite{Biggins92} and its
  ``continuous-time'' extension discussed on page~150 there.
\end{proof}

\subsection{On the largest fragment}\label{sec:A_Weierstrass}
Alike the observation made by Bertoin~\cite[Equation~(9)]{Bertoin03}
for pure homogeneous fragmentations, \cref{thm:A_Jacobi} readily reveals
the asymptotic velocity of the largest fragment~$Z_1$: if~$\PP^*$
denotes the probability~$\PP$ conditionally on non-extinction, then
\begin{flalign*}
  &&&&\lim_{t\to\infty}\frac 1t\log Z_1(t)&=\cf'(\bar q),
  &\PP^*\text{-almost surely},%
  \\[\dimexpr-\baselineskip-\lineskip\relax]
\end{flalign*}
where $\cf'(\bar q)=\cf(\bar q)/\bar q$ and provided that
$\bar q>\ubar q$.
We shall now delve deeper into the analogy with branching random walks
and tell a bit more about the asymptotic expansion of~$Z_1(t)$. To this
end, we proceed to a renormalization of the branching process~$\BRW^t$:
specifically, for
\begin{equation*}
  \BRWb^t\defeq\sum_{i=1}^\infty
  \delta_{\cf(\bar q)t-\bar q\log Z_i(t)},
\end{equation*}
which has the log-Laplace transform
\begin{equation*}
\cfb(q)\,
  \defeq\,\frac1t\log\EE\!\left[\int_{\RR}e^{-qy}\,\BRWb^t(\dd y)\right]
  =\,\cf(q\bar q)-q\cf(\bar q),\qquad q\ge0,
\end{equation*}
we are now in the so called \emph{boundary case}, namely $\cfb(0)>0$ and
$\cfb(1)=\cfb'(1)=0$. Let us also introduce the process
\begin{equation*}
  D(t)
    \defeq\int_\RR ye^{-y}\,\BRWb^t(\dd y)
    =-\bar q\left.\frac{\dd}{\dd q}\,M(t;q)\right|_{q=\bar q},
      \qquad t\ge0,
\end{equation*}
which is easily seen from the branching property to be a martingale
(rightly called the \emph{derivative martingale}) and will serve our
purpose.
\begin{corollary}\label{cor:A_Hamilton}
  Suppose~\eqref{eq:A_Dirichlet}, \eqref{eq:A_Ramanujan}, and $\bar q>\ubar q$.
  \begin{enumerate}[label=(\alph*)]
  \item Then\label{itm:Hermite}
  \begin{equation}
    \lim_{t\to\infty}\frac{\log Z_1(t)-\cf'(\bar q)t}{\log t}
      =-\frac3{2\bar q},\qquad\text{in $\PP^*$-probability}.%
    \label{eq:A_Laplace}
  \end{equation}
  \item If further~$\mu$ is non-lattice, then there exist a constant
    $C^*>0$ and a nonnegative random variable~$D_\infty$ such that, for
    every $x>0$,
  \begin{equation}
    \lim_{t\to\infty}
    \PP\!\left(t^{3/2\bar q}\,e^{-\cf'(\bar q)t}\,Z_1(t)\le x
      \right)=\,\EE\!\left[e^{-C^*D_\infty/x}\right]\!.%
    \label{eq:A_Dedekind}
  \end{equation}
    Moreover $D_\infty>0$, $\PP^*$-almost surely.\label{itm:Godel}
  \end{enumerate}
\end{corollary}
\begin{remark}
  (i)~Kyprianou et al.~\cite{Kyprianou17} recently derived an analogue
  of~\ref{itm:Hermite} for pure homogeneous fragmentations. However the
  method we employ here (for both statements) is different: basically,
  we directly transfer the known results on branching random walks to
  discrete skeletons of the growth-fragmentation, and then infer the
  behavior of the whole process with the help of \cref{lem:A_Klein} below.

  \noindent(ii)~The logarithmic fluctuations~\cite[Theorem~5.23]{Shi15}
    also show that
  \begin{equation*}
    \limsup_{t\to\infty}\frac{\log Z_1(t)-\cf'(\bar q)t}{\log t}
        \ge-\frac1{2\bar q},\qquad\text{$\PP^*$-almost surely}
  \end{equation*}
  (we conjecture that there is in fact equality), so the
  convergence~\eqref{eq:A_Laplace} cannot be strengthened.

  \noindent(iii)~Other interesting facts from the literature of
    branching random walks could be inherited. For instance, by
    specializing a recent result due to A\"id\'{e}kon and
    Shi~\cite[Theorem~5.29]{Shi15} one infers a so called
    Seneta--Heyde renormalization for the convergence of~$M(t;\bar q)$
    in \hyperref[thm:A_Jacobi2]{\cref*{thm:A_Jacobi}.(ii)}: namely
  \begin{equation*}
    \lim_{t\to\infty}\sqrt t\,M(t;\bar q)
      =\sqrt{\frac2{\pi\bar q^2\cf''(\bar q)}}\,D_\infty,%
    \qquad\text{in $\PP^*$-probability}
  \end{equation*}
  (with $D_\infty$ as above), and again this convergence cannot be
  strengthened (the $\limsup$ is infinite $\PP^*$-almost surely).
\end{remark}

\begin{lemma}[\protect{Croft--Kingman, \cite[Theorem~2]{Kingman63}}]\label{lem:A_Klein}
  Let $f\colon\oo{0}{\infty}\to\RR$ be a continuous function such that
  for every $h>0$, the sequence $f(nh),\,n\in\NN$, converges.
  Then~$f(x)$ has a limit as~$x\to\infty$.
\end{lemma}

\begin{proof}[Proof of \cref{cor:A_Hamilton}]
Let $h>0$ be any fixed time mesh. It is plain from the branching
property that~$\BRWb^{nh}$ corresponds to the individuals at generation
$n\in\NN$ of a branching random walk on~$\RR$ whose offspring point
process is distributed like~$\BRWb^h$. On the one hand, there is
\begin{equation*}
  \EE\!\left[\int_{\RR}y^2\,e^{-y}\,\BRWb^h(\dd y)\right]
  =\,h\cfb''(1)\,=\,h\bar q^2\cf''(\bar q)\,<\,\infty.
\end{equation*}
On the other hand, with the notation
$u_+\defeq\max(u,0)$ for any $u\in\RR$ and
\begin{equation*}
  X\defeq\int_\RR e^{-y}\,\BRWb^h(\dd y)=M(h;\bar q),\qquad
  \widetilde X\defeq\int_\RR y_+\,e^{-y}\,\BRWb^h(\dd y),
\end{equation*}
\cref{pro:A_Euclid} readily entails that
\begin{equation*}
  \EE\!\left[X(\log X)_+^2\right]<\,\infty\qquad\text{and}\qquad
  \EE\!\left[\widetilde X(\log\widetilde X)_+\right]<\,\infty.
\end{equation*}
(For the second, we use that
$\lvert\log(x)\rvert\le(x^\eps+x^{-\eps})/\eps$ for
every $x>0$ and any $0<\eps<\bar q-\ubar q$.)
As a result, Assumption~(H) of \cite[{\S}~5.1]{Shi15} is fulfilled%
\footnote{Strictly speaking, \cite{Shi15} also requires a
  \emph{finite} branching ($\BRWb^t(\RR)<\infty$ a.s.), but this
  condition turns out to be unnecessary (see e.g.~\cite{Mallein15};
  besides, it is not needed in the latest version of~\cite{Aidekon13}
  that we invoke to prove~\ref{itm:Godel}, and the
  conclusion of~\ref{itm:Godel} obviously implies~\ref{itm:Hermite}).}.
From Theorem~5.12 there, we obtain that for every $\eps>0$,
\begin{equation*}
  \PP\!\left(\left|\frac{\log Z_1(nh)-\cf'(\bar q)nh}{\log nh}
  +\frac3{2\bar q}\right|>\eps\right)\xrightarrow[n\to\infty]{}\,0,
  \qquad\text{for each }h>0.
\end{equation*}
As the left-hand side is a continuous function of $t\defeq nh$, the
proof of~\ref{itm:Hermite} follows from \cref{lem:A_Klein}.
Similarly, when~$\mu$ is non-lattice, then~$\BRWb^h$ is non-lattice as
well and Theorem~5.15 of~\cite{Shi15} (likewise, Theorem~1.1
of~\cite{Aidekon13}) applies: for every $y\in\RR$ and every $h>0$, the
quantity
\begin{equation*}
  \PP\!\left(\log Z_1(nh)-\cf'(\bar q)nh+\frac3{2\bar q}\log nh\le y
    \right)
\end{equation*}
has a limit as $n\to\infty$. Applying Croft--Kingman's lemma once
more gives the convergence~\eqref{eq:A_Dedekind}, where the limit is e.g.\ that
for $h=1$. In~\ref{itm:Godel}, the random variable $D_\infty$ can
be taken as the terminal value of the derivative martingale
$(D(n),\,n\in\NN)$: Theorem~5.2 of~\cite{Shi15} shows that it exists
almost surely and is positive on non-extinction. That~$D_\infty$ is
at least nonnegative holds simply because the smallest atom of~$\BRWb^n$,
\begin{equation*}
  \cf(\bar q)n-\bar q\log Z_1(n),
\end{equation*}
is bounded from below by $-\log M(n;\bar q)$, which
tends to~$\infty$ a.s.\ due to
\hyperref[thm:A_Jacobi2]{\cref*{thm:A_Jacobi}.(ii)}.
\end{proof}

\subsection{On abnormally large fragments}\label{sec:A_Descartes}
In this last section we give an estimation for the probability of
presence of fragments at scale greater than $\cf'(\bar q)$ in a
compensated fragmentation. We simply perform the very same analysis as
done in~\cite{Rouault05} for homogeneous pure fragmentations. Let us fix
two real numbers $\alpha<\beta$ and introduce
\begin{align*}
  U(t,x)
    &\defeq\PP\Bigl(\BRW^t(\cc{x+\alpha}{x+\beta})>0\Bigr),%
\\[.4em]
  V(t,x)
    &\defeq\EE\Bigl[\BRW^t(\cc{x+\alpha}{x+\beta})\Bigr],
\end{align*}
for every $t\ge0$ and $x\in\RR$.

\begin{corollary}\label{cor:A_Weyl}
  Let $q>\ubar q$. Suppose~\eqref{eq:A_Dirichlet},
  $\mu$ non-lattice, and either~\eqref{eq:A_Ramanujan} or $q\ge1$.
\begin{enumerate}[label=(\roman*)]
  \item\label{cor:A_Weyl1} Then
  \begin{equation*}
    \lim_{t\to\infty}\sqrt t\,e^{-\bigl(\cf(q)-q\cf'(q)\bigr)t}\,
    V\bigl(t,t\cf'(q)\bigr)
      =\frac{e^{-q\alpha}-e^{-q\beta}}{q\sqrt{2\pi\cf''(q)}}.
  \end{equation*}
  \item\label{cor:A_Weyl2} If further $q>\bar q$ (so that $\cf(q)-q\cf'(q)<0$),
    then
  \begin{equation*}
    \lim_{t\to\infty}\frac{U\bigl(t,t\cf'(q)\bigr)}{V\bigl(t,t\cf'(q)\bigr)}=K_q,
  \end{equation*}
    where~$K_q$ is some positive finite constant.
\end{enumerate}
\end{corollary}
\begin{remark}
  In the range $q\in\oo{\ubar q}{\bar q}$, \labelcref{cor:A_Weyl1} is the
  counterpart in mean of the convergence stated in~\cref{cor:A_Cantor} for
  $f\defeq\II_{\cc{\alpha}{\beta}}$. This convergence thus holds in
  $\Leb^1(\PP)$ thanks to the Riesz--Scheff\'{e} lemma.
\end{remark}
\begin{proof}
  The proof is a straightforward adaptation of that of Theorem~5
  in~\cite{Rouault05}. In our setting, we have
\begin{equation*}
  a\defeq\cf'(q),\qquad\Lambda^*(a)\defeq q\cf'(q)-\cf(q),
\end{equation*}
  and for any time mesh $h>0$,
\begin{equation*}
  \Lambda_h(q)\,\defeq\,\log\EE\!\left[\sum_{i=1}^\infty Z_i^q(h)\right]
    =\,h\cf(q),
\end{equation*}
  by~\eqref{eq:A_Leibniz}. From Equation~(12) in~\cite{Rouault04} we
  readily get
\begin{equation}
  \sqrt{nh}\,e^{nh\Lambda^*(a)}V(nh,anh)
    \xrightarrow[n\to\infty]{}\frac{e^{-q\alpha}-e^{-q\beta}}{q\sqrt{2\pi\cf''(q)}}.%
  \label{eq:A_Pascal1}
\end{equation}
  If furthermore $\Lambda^*(a)>0$ (i.e.\ $q>\bar q$), then
  \cref{pro:A_Euclid} ensures that the conditions of Theorem~2
  in~\cite{Rouault04} are fulfilled and therefore
\begin{equation}
  \frac{U(nh,anh)}{V(nh,anh)}\xrightarrow[n\to\infty]{} K_q^{(h)},%
  \label{eq:A_Pascal2}
\end{equation}
  where $K_q^{(h)}$ is a positive constant. Besides, the time mesh $h>0$
  in~\eqref{eq:A_Pascal1} and~\eqref{eq:A_Pascal2} is arbitrary and the left-hand
  sides are both continuous functions of the variable $t\defeq nh$.
  The existence of limits as $t\to\infty$ then comes again from
  \cref{lem:A_Klein}. (In particular, the constant~$K_q^{(h)}$
  in~\eqref{eq:A_Pascal2} does actually not depend on~$h$.)
\end{proof}

\section{Self-similar growth-fragmentations}\label{sec:A_Poincare}
As opposed to the previous part, a \emph{self-similar}
growth-fragmentation will now allow inhomogeneous fragmentation rates.
Loosely speaking, one can picture it as a homogeneous fragmentation
where each fragment is ``sped up'' all along its history by a fixed
power $\alpha\in\RR$ of its current size. If as before the Laplace
transform of the fragment sizes at genealogical births may be related
through a cumulant function~$\cf$, self-similarity induces significant
changes when we look at processes over time. Mainly, in the case
$\alpha>0$ we shall mostly focus on, where positive growth in the fragments is
thus compensated by higher dislocation rates, the typical sizes will no
longer be of exponential order (given through the derivative~$\cf'$), but
will instead encounter a polynomial decay of the type $t^{-1/\alpha}$.
Another side effect is that additive martingales appear less nicely, so
specific assumptions will be needed.

\subsection{Prerequisites}\label{sec:A_Cartan}
We begin with a quick summary of the construction and important
properties of self-similar growth-fragmentations processes. These were
introduced in~\cite{Bertoin17}; greater details as well as some
applications to random planar maps can be found in~\cite{Curien18}.

Let $\xi\defeq(\xi(t),\,t\ge0)$ be a possibly killed L\'{e}vy process and
$(\sigma^2,b,\Lambda,\mathsf{k})$ denote its characteristic quadruple
in the following sense. The Gaussian coefficient $\sigma^2\ge0$, the
drift coefficient $b\in\RR$, the L\'{e}vy measure~$\Lambda$ (that is, a
measure on~$\RR$ with $\int(1\wedge y^2)\,\Lambda(\dd y)<\infty$), and
the killing rate $\mathsf{k}\in\co{0}{\infty}$ may be recovered from
this slight variation of the L\'{e}vy--Khinchin formula:
\begin{align*}
  E\!\left[\exp\bigl(q\xi(t)\bigr)\right]=\,\exp\bigl(t\Psi(q)\bigr),\qquad\qquad%
    &\qquad t,q\ge0,%
\intertext{the Laplace exponent~$\Psi$ being written in the form}
\Psi(q)
  \defeq-\mathsf{k}+\frac12\sigma^2q^2+bq
    +\int_\RR\bigl(e^{qy}-1+q(1-e^y)\bigr)\,
    \Lambda(\dd y),&\qquad\hphantom{t,\,}q\ge0.
\end{align*}
The case $\Lambda(\oo{-\infty}{0})=0$ will be uninteresting and is
therefore excluded. We shall also assume that
$\int_{\oo{1}{\infty}}e^y\,\Lambda(\dd y)<\infty$ (which always holds
when the support of~$\Lambda$ is bounded from above), and that
\begin{equation}
  \mathsf{k}>0\qquad\text{or}\qquad
    \bigl(\mathsf{k}=0\quad\text{and}\quad\Psi'(0+)\in\co{-\infty}{0}
    \bigr)%
  \label{eq:A_Hardy}
\end{equation}
(in other words, that $\Psi(q)<0$ for some $q>0$). This latter
condition means that $\xi(t)$ either has a finite lifetime or tends to
$-\infty$ as $t\to\infty$, almost surely.
Let now $\alpha\in\RR$ and, for each $x>0$, $P_x$ be the law of the
process
\begin{gather*}
  X(t)\defeq x\exp\bigl\{\xi(\tau_{x^\alpha t})\bigr\},\qquad t\ge0,
\intertext{where}
  \tau_t\defeq\inf\left\{u\ge0\colon
    \int_0^u\exp\bigl(-\alpha\xi(s)\bigr)\,\dd s
    \ge t\right\}\!,
\end{gather*}
and with the convention that $X(t)\defeq\cem$ for
$t\ge\zeta\defeq x^{-\alpha}\int_0^\infty
  \exp(-\alpha\xi(s))\,\dd s$.
This \emph{Lamperti transform} (\cite{Lamperti72}; see also
\cite[Theorem~13.1]{Kyprianou14}) makes $(X,(P_x)_{x>0})$ be a positive
self-similar Markov process (for short, \emph{pssMp}), in the sense
that:
\begin{flalign}
&\text{For all }x>0,
&\text{the law of }\bigl(xX(x^\alpha t),t\ge0\bigr)\text{ under }P_1
 \text{ is }P_x.&&&&%
 \label{eq:A_Borel}
\end{flalign}
(Following the terminology in~\cite{Yor02}, we say
that~$X$ is a pssMp with index $1/(-\alpha)$.) Moreover,
this transformation is reversible, and since the law of the L\'{e}vy
process~$\xi$ is uniquely determined by its Laplace exponent~$\Psi$,
the pair $(\Psi,\alpha)$ characterizes the law of $X$; we call it the
characteristics of the pssMp $X$. Note that under~\eqref{eq:A_Hardy}, $X$
either is eventually absorbed to the cemetery point~$\cem$ added to the
positive half-line~$\oo{0}{\infty}$, or it converges to~$0$ as
$t\to\infty$.

The process~$X$ above will portray the typical size of a cell in the
system and is thus referred as the \emph{cell process}. Specifically, a
cell system is a process $((\CS_u,b_u),\,u\in\HUT)$ indexed on the
Ulam--Harris tree
\begin{equation*}
  \HUT\defeq\bigcup_{i=0}^\infty\NN^i,%
\end{equation*}
with the following classical notations: $\NN^0$ is reduced to the
\emph{root} of~$\HUT$, labeled~$\es$, and for any
\emph{node} $u\defeq u_1u_2\cdots u_i\in\HUT$ in this tree,
$|u|\defeq i\in\{0,1,2,\ldots\}$ refers to its \emph{generation} (or
\emph{height}), and $u1,u2,\ldots$ to its children. For each
$u\in\HUT$, $(\CS_u(t),\,t\ge0)$ is a c\`{a}dl\`{a}g process\vadjust{\goodbreak}
on
$\oo{0}{\infty}\cup\{\cem\}$ driven by~$X$ and recording the size of the
cell labeled by~$u$ since its birth time~$b_u$, which shall be
implicitly encoded in the notation~$\CS_u$. In this system $\CS_\es$
refers to Eve cell, born at time $b_\es\defeq0$, and each negative%
\footnote{In~\cite{Bertoin17}, the author only considered spectrally
negative L\'{e}vy processes so jumps were always of negative sign.
However, allowing the cells to have sudden positive growths during
their lifetimes is relevant in some applications such as those exposed
in~\cite{Curien18}. Their slightly more general setting, which we have
also chosen to adopt, does not invalidate the results
of~\cite{Bertoin17} --- the significant point being that only the
negative jumps correspond to division events while the possible
positive jumps just remain part of the trajectories of the cells.}
jump of a cell is interpreted as the birth of a daughter cell. More
precisely for every $u\in\HUT$ and $j\in\HUT$, $\CS_{uj}$ is the
process of the $j$\textsuperscript{th} daughter cell of~$u$, born at
the absolute time $b_{uj}\defeq b_u+\beta_{uj}$, where~$\beta_{uj}$ is
the instant of the $j$\textsuperscript{th} biggest positive jump%
\footnote{Recall that the cell process is either absorbed in finite time
or converges to~$0$, so the positive jumps of~$-\CS_u$ may indeed be
ranked in the decreasing order.} of~$-\CS_u$.
The law~$\Px$ of~$\CS$ is then defined recursively as the unique
probability distribution such that~$\CS_\es$ has the law~$P_x$ and,
conditionally on~$\CS_\es$, the processes
$(\CS_{iu},\,u\in\HUT),\,i\in\NN,$ are independent with respective laws
$\Px[x_i],\,i\in\NN$, where $(x_1,\beta_1),(x_2,\beta_2),\ldots$ is the
sequence of positive jump sizes and times of $-\CS_\es$, sorted by
decreasing sizes (with $\beta_i<\beta_{i+1}$ if $x_i=x_{i+1}$). Here, we
agree that $x_i\defeq\cem$ and $\beta_i\defeq\infty$ if~$\CS_\es$ has
less than~$i$ negative jumps, and we let~$\Px[\cem]$ denote the law of
the degenerate cell system where $\CS_u\equiveq\cem$ for every
$u\in\HUT$, $b_\es\defeq0$ and $b_u\defeq\infty$ for $u\neq\es$.

The associated growth-fragmentation process is the process of the
family of (the sizes of) all alive cells in the system:
\begin{align*}
&&&&&&&&&&&&\bs X(t)&\defeq\Bigl\{\!\!\Bigl\{\CS_u(t-b_u)\colon
    u\in\HUT,\,b_u\le t<d_u\Bigr\}\!\!\Bigr\},&\qquad t\ge0\hphantom{.}&&
\intertext{(with~$b_u$ and~$d_u$ denoting respectively the birth time and the death
time of the cell labeled by $u$).
Additionally to the scaling parameter $\alpha$, one other specific
quantity is}
&&&&&&&&&&&&\cf(q)&\defeq\Psi(q)+\int_{\oo{-\infty}{0}}(1-e^y)^q\,\Lambda(\dd y),%
    &\qquad q\ge0.
\end{align*}
If $\alpha=0$ and~$\Lambda$ has support in
$\cc{-\log2}{0}$, then \cite[Proposition~3]{Bertoin17}
$\bs X$ is merely a compensated fragmentation of the type considered
in \cref{sec:A_Riemann}, and the notation~$\cf$ there is compatible with the
one we use here: more precisely, $\bs X$ has diffusion
coefficient~$\sigma^2$, growth rate~$b$ and dislocation measure
$\nu\defeq\textsf{k}\,\delta_{\bs 0}+\nu_2$, where
$\bs0\defeq(0,0,\ldots)\in\MP$ is the null mass-partition and~$\nu_2$ is
the image of~$\Lambda$ by the map
$x\in\cc{-\log 2}{0}\mapsto(e^x,1-e^x,0,\ldots)\in\MP$ (the
fragmentation is binary).

We shall work under the assumption
\begin{equation}
\exists q\ge0,\ \cf(q)\le0,\label{eq:A_Kolmogorov}
\end{equation}
see \cite[Theorem~2]{Bertoin17}. Then for each
time~$t$, the family~$\bs X(t)$ may be ranked in the non-increasing
order, i.e.\ $\bs X(t)\defeq(X_1(t),X_2(t),\ldots)$ with
$X_1(t)\ge X_2(t)\ge\cdots\ge0$. Further, the self-similarity
property~\eqref{eq:A_Borel} extends to the process
$\bs X\defeq(\bs X(t),\,t\ge0)$, and there is the branching property.
Formally, if~$\Pb$ denotes the law of~$\bs X$ under~$\Px$, then firstly,
for every $x>0$, the distribution of $(x\bs X(x^\alpha t),\,t\ge0)$
under~$\Pb[1]$ is~$\Pb$, and secondly, for each $s\ge0$, conditionally
on $\bs X(s)=(x_1,x_2,\ldots)$, the process $(\bs X(t+s),\,t\ge0)$ is
independent of $(\bs X(r),\,0\le r\le s)$ and has the same law
as the non-increasing rearrangement of the family
$(X^{(i)}_j,\,i,j\in\NN)$, where the~$\bs X^{(i)}$ are independent
self-similar growth-fragmentations with respective laws~$\Pb[x_i]$.

In the sequel we mainly focus on large time asymptotics for the
growth-fragmentation process~$\bs X$. Since we can refer to \cref{sec:A_Riemann} when
$\alpha=0$, and because the growth-fragmentation is eventually extinct
when the scaling parameter~$\alpha$ is negative
\cite[Corollary~3]{Bertoin17}, we will mostly suppose $\alpha>0$. Note
in this case that~\eqref{eq:A_Kolmogorov} is a necessary and sufficient
condition preventing local explosion of the
fragmentation~\cite{Stephenson16}, that is a phenomenon causing
infinitely many particles of arbitrary large sizes to be produced in
almost surely finite time (which in particular would impede us to list
the elements of~$\bs X(t)$ in the non-increasing order). Like in
\cref{sec:A_Riemann}, the function
$\cf\colon\co{0}{\infty}\to\oc{-\infty}{\infty}$ will be of greatest
importance in the study. It is clearly convex; therefore the equation
$\cf(q)=0$ has at most two solutions. We assume from here on that these
two solutions exist --- more precisely that the
\emph{Malthusian hypotheses} hold:
\begin{equation}
  \text{there exist }0<\cffr<\cfsr
  \text{ such that }\quad
    \cf(\cffr)=\cf(\cfsr)=0\ \text{ and }\
    \cf'(\cffr)>-\infty%
  \label{eq:A_Liouville}
\end{equation}
(note then that $\cf'(\cffr)<0$, by convexity).
Condition~\eqref{eq:A_Liouville} implies that $\cf(q)<0$ for some $q>0$, which
in turn implies~\eqref{eq:A_Kolmogorov}, and~\eqref{eq:A_Hardy} (because $\Psi\le\cf$).

As before, limit theorems for the growth-fragmentation process
$\bs X$ will involve the terminal value of some additive martingale,
namely the \emph{Malthusian martingale}
\begin{equation*}
  M^-(t)\defeq\sum_{i=1}^\infty X_i^\cffr(t),\qquad t\ge0.
\end{equation*}
In this direction, results of~\cite{Curien18} will be of fundamental
use; we restate some of them here for sake of reference.
\begin{proposition}[from \protect{\cite[Theorem~3.10.(ii),
  Corollaries~3.7.(ii) and 3.9]{Curien18}}]\label{pro:A_Eisenstein}\leavevmode

  \noindent Suppose $\alpha>0$.
\begin{enumerate}[label=(\roman*)]
    \item The process $(M^-(t),\,t\ge0)$ under~$\Pb$ is a uniformly
      integrable martingale; more precisely it is bounded in
      $\Leb^p(\Pb)$ for every $1<p<\cfsr/\cffr$.%
    \label{pro:A_Eisenstein1}
    \item For every $0<q<(\cfsr-\cffr)/\alpha$, the process
    \begin{equation*}
      \sum_{i=1}^\infty X_i^{q\alpha+\cffr}(t),\qquad t\ge0,
    \end{equation*}
       is a supermartingale converging to~$0$ in $\Leb^1(\Pb)$: more
       precisely,
    \begin{equation*}
      \Eb\!\left[\sum_{i=1}^\infty
      X_i^{q\alpha+\cffr}(t)\right]\sim\,c(q)\,x^\cffr\,t^{-q}
    \end{equation*}
      as $t\to\infty$, for some constant $c(q)>0$.%
     \label{pro:A_Eisenstein2}
\end{enumerate}
\end{proposition}
\begin{remark}\label{rem:A_Bernoulli}
  We find relevant to mention that~\cite{Curien18} also introduced the
  genealogical martingale
\begin{flalign*}
  &&&&&&&&\mathcal{M}^-(n)&\defeq\!\!\sum_{|u|=n+1}\!\!\!\CS_u^\cffr(0),%
    &n\ge0,&&&&&&
\intertext{called the \emph{intrinsic martingale}, which under~$\Px$ is always
  uniformly integrable. When $\alpha\ge0$, there is the remarkable fact}
  &&&&&&&&M^-(t)&\hphantom{:}=\hphantom{:}
    \Ex\!\left[\mathcal{M}^-(\infty)\;\middle|\;\mathcal F_t\right]\!,%
    &t\ge0,
\end{flalign*}
  with $(\mathcal F_t)_{t\ge0}$ the canonical filtration of~$\bs X$.
  In particular, $M^-(\infty)=\mathcal{M}^-(\infty)$ almost surely.
\end{remark}

Additive martingales --- in the present context, the Malthusian
martingale $(M^-(t))_{t\ge0}$, are of important interest since the
celebrated work of Lyons et al.~\cite{Lyons95}. Roughly speaking,
one can perform a change of probability measure in terms of
the terminal value $M^-(\infty)$ so that the genealogical system may be
observed from the point of view of a randomly tagged branch.
Specifically, write $\partial\HUT$ for the set of \emph{leaves}
of~$\HUT$, each of which determines a unique branch from the root. For
every leaf $\ell\in\partial\HUT$, let $\ell(n)$ denote its unique
ancestor at generation $n\ge0$, and
$\CS_\ell\defeq(\CS_\ell(t),\,t\ge0)$ be the process of the cell
on the branch from~$\es$ to~$\ell$:
\begin{equation*}
  \CS_\ell(t)\defeq\CS_{\ell[t]}(t-b_{\ell[t]}),\qquad t\ge0,
\end{equation*}
where $\ell[t]$ labels the cell in this branch which is alive at
time~$t$ (i.e.\ $\ell[t]$ is the unique ancestor~$u$ of~$\ell$ such that
$b_{u}\le t<b_{\ell(|u|+1)}$), with the convention that
$\CS_\ell(t)\defeq\partial$ for
$t>\mathop{\lim\!\!\uparrow}_{n\to\infty}b_{\ell(n)}\eqdef b_\ell$. Next we consider a
random leaf $\RL\in\partial\HUT$ and we define for every $x>0$ the joint
distribution~$\Pm$ of $(\CS,\RL)$ as follows. Under~$\Pm$, the law of
$\CS\defeq(\CS_u,\,u\in\HUT)$ is absolutely continuous with respect
to~$\Px$ with density $x^{-\cffr}M^-(\infty)$, and the law of~$\RL$
conditionally on~$\CS$ is
\begin{equation}
  \Pm(u\text{ ancestor of }\RL\mid\CS)
  \defeq\lim_{n\to\infty}\frac1{\mathcal{M}^-(\infty)}
    \sum_{|v|=n}\CS_{uv}^\cffr(0).%
  \label{eq:A_Kepler}
\end{equation}
Denoting $\TC\defeq\CS_\RL$ the \emph{randomly tagged cell}, Bertoin
et al.~\cite{Curien18} derived:
\begin{proposition}[\protect{from \cite[Theorem~4.7 and
  Proposition~4.6]{Curien18}}]\leavevmode
\begin{enumerate}[label=(\roman*)]
  \item The process $(\TC,(\Pm)_{x>0})$ is a pssMp with characteristics
      $(\Phi^-,\alpha)$, where
  \begin{equation}
    \Phi^-(q)\defeq\cf(q+\cffr),\qquad q\ge0.\label{eq:A_Banach}
  \end{equation}
  \item \emph{Many-to-one formula.}
    For every $x>0$, every $t\ge0$, and every measurable function
    $f\colon\oo{0}{\infty}\to\oo{0}{\infty}$, we have
  \begin{equation}
    \Eb\!\left[\sum_{i=1}^\infty X_i^\cffr(t)\,f\bigl(X_i(t)\bigr)\right]
      =\,x^\cffr\,\Em\!\left[f\bigl(\TC(t)\bigr)\right]\!,%
    \label{eq:A_Lebesgue}
  \end{equation}
    with the convention $f(\cem)\defeq0$.
\end{enumerate}
\end{proposition}
Formula~\eqref{eq:A_Lebesgue} will be a key ingredient for our purpose.
Roughly speaking, it says that the intensity of the weighted point
measure $\sum_{z\in\bs X}z^\cffr\delta_z$ is captured by the law of the
randomly tagged cell~$\TC$ (hence the denomination ``many-to-one'').

When $\alpha>0$, unlike in the homogeneous case, a \emph{polynomial}
decrease in the size of the fragments is expected. Large-time
asymptotics for their empirical measure will be retrieved in the next
section. In \cref{sec:A_Hadamard} we center our attention on the largest
fragment. Lastly, in \cref{sec:A_Peano}, we discuss the convergence of
the empirical measure of the fragments taken at the instant when they
become smaller than a vanishing threshold.

\subsection{Convergence of the empirical measure}
We are here especially interested in the convergence of the empirical
measure~$\rho_t^{(\alpha)}$ given by
\begin{equation*}
  \bigl\langle\rho_t^{(\alpha)},f\bigr\rangle\defeq\sum_{i=1}^\infty
    X_i^\cffr(t)\,f\bigl(t^{1/\alpha}X_i(t)\bigr),
\end{equation*}
for $\alpha>0$. From now on, we shall suppose that the L\'{e}vy
process~$\tc$ associated with the tagged cell~$\TC$ via Lamperti's
transformation is not \emph{arithmetic}, in the sense that there is no
$r>0$ such that $\PP(\tc(t)\in r\ZZ)=1$ for all $t\ge0$.
To state our result, let us define the probability distribution~$\rho$
on~$\oo{0}{\infty}$ by
\begin{equation*}
\int_0^\infty f(y)\,\rho(\dd y)\defeq
  \frac{-1}{\alpha\cf'(\cffr)}\,E\!\left[I^{-1}\,f\bigl(I^{1/\alpha}
    \bigr)\right]\!,
\end{equation*}
where
\begin{equation}
  I\defeq\int_0^\infty\exp\bigl(\alpha\tc(s)\bigr)\,\dd s%
  \label{eq:A_Weil}
\end{equation}
is the so called \emph{exponential functional} of $\alpha\tc$.
The following completes the results of Bertoin~\cite{Bertoin03}
and Bertoin and Gnedin~\cite{Gnedin04} relative to self-similar pure
fragmentations, and differs substantially from the homogeneous case
(\cref{cor:A_Cantor}).
\begin{theorem}\label{thm:A_Monge}
  For every $1<p<\cfsr/\cffr$ and for every bounded continuous function
  $f\colon\oo{0}{\infty}\to\RR$,
\begin{equation*}
  \lim_{t\to\infty}
  \bigl\langle\rho_t^{(\alpha)},f\bigr\rangle
    =M^-(\infty)\int_0^\infty f(y)\,\rho(\dd y),%
  \qquad\text{in }\Leb^p(\Pb[1]).
\end{equation*}
  Consequently, the random measure $\rho_t^{(\alpha)}$ converges in
  $\Pb[1]$-probability to $M^-(\infty)\,\rho$ as $t\to\infty$, in the
  space of finite measures on~$\oo{0}{\infty}$ endowed with the
  topology of weak convergence.
\end{theorem}
\begin{remark}
  Note the presence of the random factor $M^-(\infty)$, which does not
  appear in~\cite{Bertoin03} because the Malthusian martingale is
  trivial for \emph{conservative} pure fragmentations%
  \footnote{Since then the cumulant function vanishes at~$1$ and the
  total mass of the fragments at any generation is constant.}.
  It does nonetheless appear in the non-conservative
  case~\cite{Gnedin04}; however the method used there leads to a
  $\Leb^2$-convergence that we cannot hope for growth-fragmentations
  when $\cfsr/\cffr<2$, and it seems anyway difficult to generalize.
\end{remark}

Exponential functionals of L\'{e}vy processes such as~\eqref{eq:A_Weil} arise
in a variety of contexts and their laws have been widely
studied, see the survey~\cite{Yor05} and the recent
works~\cite{Pardo12,Patie12,Pardo13,Arista15}. In particular, Pardo
et al.~\cite{Pardo12} showed that under mild assumptions, they can be
factorized into the product of two independent exponential functionals
associated with companion L\'{e}vy processes, and the distributions of both
these functionals are uniquely determined by either their positive or
their negative moments. To name just one example, in the common
situation where~$\tc$ is spectrally negative
($\Lambda(\oo{0}{\infty})=0$), we have
\cite[Corollary~2.1]{Pardo12} that $I\eqlaw J/\Gamma$, with~$J$ the
exponential functional of the descending ladder height process
of~$\alpha\tc$ and~$\Gamma$ an independent Gamma random variable with
parameter $(\cfsr-\cffr)/\alpha$. Further, the density of~$I$ has a
polynomial tail of order $1+(\cfsr-\cffr)/\alpha$ and admits a
semi-explicit series expansion.

The distribution of~$I$ (likewise, $\rho$) naturally takes part
in asymptotics of the tagged cell~$\TC$:
\begin{lemma}\label{lem:A_Grassmann}
As $t\to\infty$, the random variable $t^{1/\alpha}\TC(t)$ under~$\Pm[1]$
converges in distribution to~$\rho$. Moreover,
$\int_0^\infty y^{q\alpha}\,\rho(\dd y)<\infty$ for every
$0\le q<1+(\cfsr-\cffr)/\alpha$.
\end{lemma}
\begin{proof}
Clearly, $(1/\TC(t),\,t\ge0)$ is a pssMp with self-similarity
index~$1/\alpha$ associated with~$-\tc$, where~$\tc$ is a L\'{e}vy process
with the Laplace exponent~$\Phi^-$ in~\eqref{eq:A_Banach}. According
to~\cite[Theorem~1]{Yor02}, all we need to check to prove the
first part of the statement is that~$-\tc(1)$ admits a finite and
positive first moment, which is implied by the Malthusian
hypotheses~\eqref{eq:A_Liouville}: indeed,
\begin{equation*}
  E[-\tc(1)]=-(\Phi^-)'(0+)=-\cf'(\cffr)\in\oo{0}{\infty}.
\end{equation*}
The existence of moments is quite straightforwardly adapted
from the proof of~\cite[Theorem~3]{Yor05}.
\end{proof}

We now turn to the proof of \cref{thm:A_Monge}, arguing along the lines of
\cite[Theorem~1.3]{Bertoin06}. The main idea is that, by the
branching and scaling properties, the empirical measure of the fragments
can be rewritten as the sum of identically distributed pieces arising
from an intermediate (arbitrary large) time, which are all
independent conditionally on the past. With the help of a (conditional)
law of large numbers, we are then reduced to a first moment estimate for
some additive functional of the growth-fragmentation, which we can work
out thanks to the many-to-one formula and the asymptotic behavior of the
tagged fragment above.

\begin{proof}[Proof of \cref{thm:A_Monge}]
  Using the branching property at time~$t$ and the self-similarity
  of~$\bs X$ we can write, on the event $\{\bs X(t)=(x_1,x_2,\ldots)\}$,
  \begin{align}
    \bigl\langle\rho^{(\alpha)}_{t+t^2},f\bigr\rangle
    &\defeq\sum_{i=1}^\infty X_i^\cffr\bigl(t+t^2\bigr)\,f\!
      \left(\bigl(t+t^2\bigr)^{1/\alpha}X_i\bigl(t+t^2\bigr)\right)
    =\,\sum_{i=1}^\infty\lambda_i(t)Y_i(t),\notag
  \intertext{where $\lambda_i(t)\defeq X_i^\cffr(t)=x_i^\cffr$ and}
    Y_i(t)&\defeq\sum_{j=1}^\infty X_{i,j}^\cffr
     \bigl(x_i^\alpha\,t^2\bigr)\,f\!\left(\bigl(t+t^2\bigr)^{1/\alpha}
       x_iX_{i,j}\bigl(x_i^\alpha\,t^2\bigr)\right)\!,%
    \label{Poisson}
  \end{align}
  the families $(X_{i,1},X_{i,2},\ldots),\,i\ge1,$ being
  i.i.d.\ copies independent of~$\bs X$, having all the same law~$\Pb[1]$.
  Clearly, the~$Y_i$ are independent conditionally on
  $\bs\lambda(t)\defeq(\lambda_1(t),\lambda_2(t),\ldots)$ and, if we
  introduce
\begin{equation*}
  \bar Y_i\defeq\|f\|_\infty\,\sup_{t\ge0}\sum_{j=1}^\infty
    X_{i,j}^\cffr(t),
\end{equation*}
  then thanks to \hyperref[pro:A_Eisenstein1]{\cref*{pro:A_Eisenstein}.(i)} and
  Doob's maximal inequality, the~$\bar Y_i$ are i.i.d.\ random variables
  in $\Leb^p(\Pb[1])$ such that $|Y_i(t)|\le\bar Y_i$ for
  all $t\ge0$. For the same reason,
\begin{align*}
  \sup_{t\ge0}\Eb[1]\!\left[\left(\sum_{i=1}^\infty
  \lambda_i(t)\right)^{\!p\,}\right]&<\,\infty
\intertext{and further, using \hyperref[pro:A_Eisenstein2]{\cref*{pro:A_Eisenstein}.(ii)},}
  \lim_{t\to\infty}\Eb[1]\!\left[\sum_{i=1}^\infty\lambda_i^p(t)\right]
    &=\,0.
\end{align*}
  By a variation of the law of large numbers (\cite{Nerman81}; see also
  \cite[Lemma~1.5]{Bertoin06}) we then have
\begin{align}
&&&&\notag\lim_{t\to\infty}\sum_{i=1}^\infty\lambda_i(t)
  \Bigl(Y_i(t)-\Eb[1][Y_i(t)\mid\bs\lambda(t)]\Bigr)=0,
&\qquad\text{in }\Leb^p(\Pb[1]).&&
\intertext{Consequently, the proof boils down to showing that}
&&&&\lim_{t\to\infty}\sum_{i=1}^\infty\lambda_i(t)\,
      \Eb[1][Y_i(t)\mid\bs\lambda(t)]
    =M^-(\infty)\int_0^\infty f(y)\,\rho(\dd y),&\qquad%
\text{in }\Leb^p(\Pb[1]),
\label{eq:A_Jordan}
\end{align}
  where, applying the many-to-one formula~\eqref{eq:A_Lebesgue},
\begin{equation*}
  \Eb[1][Y_i(t)\;\big|\;\bs\lambda(t)]
  =\Em[1]\!\left[f\!\left((1+t^{-1})%
      ^{1/\alpha}\,x_i\,t^{2/\alpha}\,
        \TC(x_i^{\alpha}t^2)\right)\right]\!.
\end{equation*}
  But we know from \cref{lem:A_Grassmann} that as $s\to\infty$, the law of
  $s^{1/\alpha}\TC(s)$ under~$\Pm[1]$ converges weakly to~$\rho$.
  On the one hand, it thus follows that
  \begin{equation*}
    \Em[1]\!\left[f\!\left((1+t^{-1})%
      ^{1/\alpha}\,x_i\,t^{2/\alpha}\,
        \TC(x_i^\alpha t^2)\right)\right]\xrightarrow[t\to\infty]{}%
    \int_0^\infty f(y)\,\rho(\dd y)
  \end{equation*}
  uniformly in~$i$ such that, say, $x_i^\alpha t^2>\sqrt t$, i.e.\
  $x_i>t^{-3/2\alpha}$. On the other hand, applying
  again~\eqref{eq:A_Lebesgue}, the quantity
  \begin{equation}
    \sum_{i=1}^\infty X_i^\cffr(t)\,\II_{\{X_i(t)\le t^{-3/2\alpha}\}}
    \label{eq:A_Fourier}
  \end{equation}
  has, under~$\Pb[1]$, mean
  \begin{equation*}
    \Pm[1]\!\left(t^{1/\alpha}\TC(t)<t^{-1/2\alpha}\right)\!,
  \end{equation*}
  which tends to~$0$ as $t\to\infty$. Since~\eqref{eq:A_Fourier} is bounded
  in $\Leb^q(\Pb[1])$ for every $p<q<\cfsr/\cffr$, it also converges
  to~$0$ in $\Leb^p(\Pb[1])$ (by H\"{o}lder's inequality). Putting
  everything together yields~\eqref{eq:A_Jordan}, and thus the first part of
  the statement.

  The second part is derived from standard arguments: the space
  $\Co_c(\oo{0}{\infty})$ of continuous functions on~$\oo{0}{\infty}$
  with compact support being separable, a diagonal extraction procedure
  easily entails, for every sequence $t_n\to\infty$, that there exists
  an extraction $\sigma\colon\NN\to\NN$ such that, almost surely,
  \begin{equation*}
    \forall f\in\Co_c(\oo{0}{\infty}),\quad
      \bigl\langle\rho_{t_{\sigma(n)}}^{(\alpha)},f\bigr\rangle
        \xrightarrow[n\to\infty]{} M^-(\infty)\int_0^\infty f(y)\,\rho(\dd y),
  \end{equation*}
  i.e.\ $\rho_{t_{\sigma(n)}}^{(\alpha)}$ converges vaguely to
  $M^-(\infty)\,\rho$, a.s. Since the total mass is conserved, that is
  \begin{equation*}
    \bigl\langle\rho_t^{(\alpha)},1\bigr\rangle
      =\sum_{i=1}^\infty X_i^\cffr(t)
    \xrightarrow[t\to\infty]{}M^-(\infty)
      =\bigl\langle M^-(\infty)\,\rho,1\bigr\rangle\qquad\text{a.s.},
  \end{equation*}
  the convergence of $\rho_{t_{\sigma(n)}}^{(\alpha)}$ toward
  $M^-(\infty)\,\rho$ is actually weak. The conclusion follows easily.
\end{proof}

The existence of moments for~$\rho$ (\cref{lem:A_Grassmann}) allows us
to strengthen \cref{thm:A_Monge}:
\begin{corollary}
  For every $0<q<(\cfsr-\cffr)/\alpha$, every measurable function
  $f\colon\oo{0}{\infty}\to\RR$ such that $f(y)=O(y^{q\alpha})$, and
  every $1<p<\cfsr/(q\alpha+\cffr)$,
\begin{equation*}
  \lim_{t\to\infty}
    \bigl\langle\rho_t^{(\alpha)},f\bigr\rangle
  =M^-(\infty)\int_0^\infty f(y)\,\rho(\dd y),%
\qquad\text{in }\Leb^p(\Pb[1]).
\end{equation*}
\end{corollary}
\begin{proof}
  Approximating $y\mapsto y^{-q\alpha}f(y)$ by simple functions, it is
  enough to do the proof for $f=f_q\colon y\mapsto y^{q\alpha}$, that
  is to prove:
\begin{equation*}
  \lim_{t\to\infty}t^q\sum_{i=1}^\infty X_i^{q\alpha+\cffr}(t)
    =M^-(\infty)\int_0^\infty y^{q\alpha}\,\rho(\dd y),%
  \qquad\text{in }\Leb^p(\Pb[1]).
\end{equation*}
  This is of course not a direct consequence to \cref{thm:A_Monge}
  because~$f_q$ is not a bounded continuous function; nevertheless we
  can repeat the argument used in the previous proof. Observing that
  $q\alpha+\cffr\in\oo{\cffr}{\cfsr}$ and defining $Y_i(t)$ as
  in~\eqref{Poisson} but with~$f_q$ in place of~$f$, we easily check with the
  help of~\cref{pro:A_Eisenstein} and H\"{o}lder's inequality that conditionally
  on~$\bs X$, the~$Y_i$ are independent supermartingales bounded in
  $\Leb^p(\Pb[1])$ for every $1<p<\cfsr/(q\alpha+\cffr)$. Therefore, all
  that remains to show is the convergence
\begin{equation*}
  t^q\,\Eb[1]\!\left[\sum_{i=1}^\infty X_i^{q\alpha+\cffr}(t)\right]
    =\,\Em[1]\!\left[\left(t^{1/\alpha}\TC(t)\right)^{\!q\alpha\,}
      \right]
    \xrightarrow[t\to\infty]{}\int_0^\infty y^{q\alpha}\,\rho(\dd y),
\end{equation*}
  where the equality is just an application of the many-to-one
  formula~\eqref{eq:A_Lebesgue}. Since we already know that
  $t^{1/\alpha}\TC(t)$ converges in distribution toward~$\rho$
  (\cref{lem:A_Grassmann}), it suffices to show that
  $((t^{1/\alpha}\TC(t))^{q\alpha})_{t\ge0}$ is bounded in
  $\Leb^r(\Pm[1])$ for some $r>1$, which is immediate using again the
  many-to-one formula and the convergence rate in
  \hyperref[pro:A_Eisenstein2]{\cref*{pro:A_Eisenstein}.(ii)}
  (we can take $1<r<(\cfsr-\cffr)/q\alpha$).
\end{proof}

\subsection{Asymptotic behavior of the largest fragment}\label{sec:A_Hadamard}
For pure self-similar fragmentations with scaling parameter $\alpha>0$,
it is known~\cite{Bertoin03} that the size of the largest fragment
decreases like $t^{-1/\alpha}$ as $t\to\infty$. The same holds for
growth-fragmentations%
\footnote{For simplicity, we suppose that the cell
process has no positive jumps, though this restriction is probably
superfluous.}:
\begin{theorem}\label{thm:A_Sznitman}
  Assume again~\eqref{eq:A_Liouville}, $\alpha>0$, and that~$\tc$ is not
  arithmetic, and suppose further that $\Lambda(\oo{0}{\infty})=0$.
  Let $S\defeq\{\forall t\ge0,\ \bs X(t)\neq\es\}$ be the non-extinction
  event, and $\Px[]^*\defeq\Px[1](\,\cdot\mid S)$. Then
\begin{equation*}
  \lim_{t\to\infty}\frac{\log X_1(t)}{\log t}
    =-\frac1{\alpha},\qquad\text{in $\Px[]^*$-probability}.
\end{equation*}
\end{theorem}

\begin{figure}
  \centering\includegraphics[width=\textwidth]{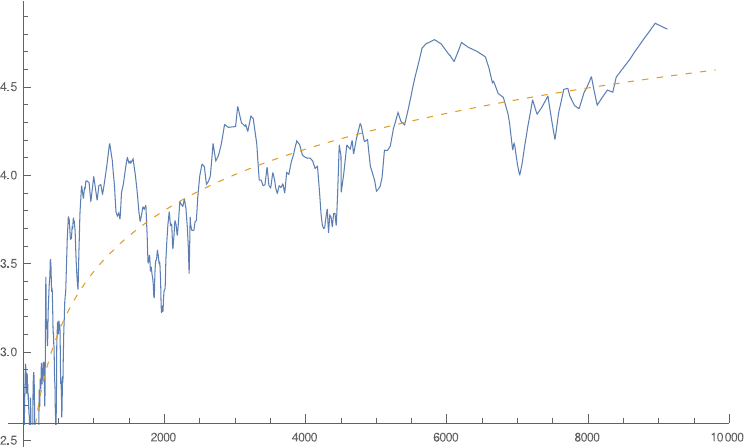}%
  \caption{Simulation of $-\log X_1$ in a self-similar
    growth-fragmentation process with scaling $\alpha=2$.
  (The dashed line represents the map $t\mapsto\frac1{\alpha}\log t$.)}
\end{figure}

\begin{proof}[Proof of the lower bound]
The fact that the $\Px[]^*$-$\liminf$ of $\log X_1(t)/\log t$ as
$t\to\infty$ is at least~$-1/\alpha$ follows by comparison
with the randomly tagged cell~$\TC$. Indeed, we know by \cref{lem:A_Grassmann}
that $\log\TC(t)/\log t$ converges to $-1/\alpha$ in
$\Pm[1]$-probability. Because~$X_1(t)$ is the size of the largest
fragment and~$\TC(t)$ is that of some other fragment in the system, we
deduce that for every $\eta>0$,
\begin{equation*}
  \Pm[1]\!\left(\frac{\log X_1(t)}{\log t}
    +\frac1{\alpha}<-\eta\right)
  \le\,\Pm[1]\!\left(\frac{\log \TC(t)}{\log t}
    +\frac1\alpha<-\eta\right)
  \xrightarrow[t\to\infty]{}\,0.
\end{equation*}
Since~$\dd\Pm[1]/\dd\Px[1]=M^-(\infty)$, which by the branching
property is positive $\Px[1]$-a.s.\ on~$S$, the latter convergence also
holds with~$\Px[]^*$ in place of~$\Pm[1]$.
\end{proof}

For the other direction, we need to make sure that asymptotically,
if the largest fragment ever exceeds the level $t^{-1/\alpha}$, it is
unlikely that one of its parents has gone far below this level
before~$t$. To this end, we write
$X_1(t)\defeq\CS_{u^*(t)}(t-b_{u^*(t)})$ with
$u^*(t)\defeq\argmax_{u\in\HUT,b_u\le t<d_u}\CS_u(t-b_u)$ (in case of
\emph{ex aequo}, we choose $u^*(t)$ to be minimal in lexicographic
order), and introduce the event
\begin{equation*}
  \mathcal H_t(\eps)\colon\quad
  \CS_v(s-b_v)<\eps,\enspace
    \text{ for some time~$s$ and ancestor $v\preccurlyeq u^*(t)$ with }
      b_v\le s<d_v\wedge t.
\end{equation*}
The following statement is tailored for our purpose.
\begin{proposition}\label{pro:A_Lie}
  There exists $\theta\ge\cfsr$ such that
  \begin{gather}
    \left\{\,\mathop{\sup}_{x>1}\;\;,\,
    \mathop{\vphantom{\sup}\lim}_{x\to\vphantom{1}\infty}\,\right\}\,
        \frac1{\log x}\log\Px[1]\!\left(\sup_{s\ge0}X_1(s)>x\right)
      =\,-\theta.%
    \label{eq:A_Einstein}
  \intertext{Furthermore, for every $\gamma,\eps\in\oo{0}{1}$
    and every $t>0$,}
    \Px[1]\Bigl(\mathcal H_t(\eps),\,
      X_1(t)>\eps^\gamma\Bigr)\,
      \le\,\eps^{(1-\gamma)\theta}.\notag
  \end{gather}
\end{proposition}
\begin{proof}
  We may assume $\alpha=0$ as the statement does not depend
  on~$\alpha$. The first
  assertion is a large deviation estimate for the probability~$F(x)$
  that $T^+(x)<\infty$, where
  $T^+(x)\defeq\inf\{s\ge0\colon X_1(s)>x\}$. To eventually obtain a
  fragment larger than~$xy$ in the growth-fragmentation, for $x,y>1$,
  it is enough that the largest particle~$X_1$ first reaches some
  level~$z>x$, and that the subsequent fragmentation of this particle
  produces a fragment with size larger than~$xy$. But by scaling we
  have, for any $z>x$,
\begin{equation*}
    \Px[z]\!\left(T^+(xy)<\infty\right)
      =\,\Px[1]\!\left(T^+\!\left(xy/z\right)<\infty\right)
      =\,F\bigl(xy/z\bigr)\,
      \ge\,F(y)=\,\Px[1]\!\left(T^+(y)<\infty\right)\!,
\end{equation*}
  so that the branching property at
  $T^+(x)$ yields, since $z\defeq X_1(T^+(x))>x$ on the event
  $\{T^+(x)<\infty\}$,
\begin{equation*}
  F(xy)
    \,=\,\Px[1](T^+(xy)<\infty)\,
    \ge\,\Ex[1]\!\left[\II_{\{T^+(x)<\infty\}}
      F\!\left(xy/X_1\bigl(T^+(x)\bigr)\right)\right]\ge F(x)F(y).
\end{equation*}
  \cref{eq:A_Einstein} then arises from the subadditive lemma
  (see e.g.\ \cite[Theorem~16.2.9]{Kuczma09}). The lower bound
  $\theta\ge\cfsr$ is just a consequence of Doob's maximal inequality
  applied to the process
\begin{equation*}
    M^+(s)\defeq\sum_{i=1}^\infty X_i^\cfsr(s),\qquad s\ge0,
\end{equation*}
  which \cite[Corollary~3.7.(i)]{Curien18} is a martingale
  (for $\alpha=0$): namely
\begin{equation*}
  F(x)=\Px[1]\!\left(\sup_{s\ge0}X_1(s)>x\right)
    \le\,\Px[1]\!\left(\sup_{s\ge0}M^+(s)>x^{\cfsr}\right)
    \le\,x^{-\cfsr},\qquad x>1.
\end{equation*}
  Next, we take $x\defeq\eps^{\gamma-1}$
  and apply again the scaling property: we deduce that, for every
  $0<y<\eps$,
  \[F(\eps^\gamma/y)\,
    =\,\Px[y]\!\left(\sup_{s\ge0} X_1(s)>\eps^\gamma\right)
    \le\,\eps^{(1-\gamma)\theta}.\]
  But the event $\mathcal{H}_t(\eps)$ holds precisely when the
  cell process following the ancestral lineage of~$u^*(t)$ has reached
  a value $0<y<\eps$ before~$t$. Using the branching property at
  the first time this happens, the second assertion is then easily
  proved.
\end{proof}
We can now derive the upper bound and complete the proof
of~\cref{thm:A_Sznitman}.
\begin{proof}[Proof of the upper bound]
  Let $0<\eta<1$ and observe that
  $\delta\defeq\eta-(1-\eta)(1-\gamma)/\gamma$ lies in~$\oo{0}{\eta}$
  for any $\gamma\in\oo{1-\eta}{1}$ arbitrarily fixed. Define
  $\eps\defeq t^{-(1-\delta)/\alpha}$ for $t>1$, so that
  $\eps^\gamma=t^{-(1-\eta)/\alpha}$, and
  \begin{align*}
    \Px[1]\!\left(X_1(t)>t^{-(1-\eta)/\alpha}\right)
      &=\,\Px[1]\!\left(X_1(t)>\eps^\gamma\right)\\
      &=\,\Px[1]\Bigl(\mathcal H_t(\eps),\,
      X_1(t)>\eps^\gamma\Bigr)
      \,+\,\Px[1]\Bigl(\mathcal H_t(\eps)^\complement,\,
      X_1(t)>\eps^\gamma\Bigr).
  \end{align*}
  By \cref{pro:A_Lie},
  \begin{equation*}
    \Px[1]\Bigl(\mathcal H_t(\eps),\,
      X_1(t)>\eps^\gamma\Bigr)\,
      \le\,\,t^{-(1-\delta)(1-\gamma)\theta/\alpha}
      \,\xrightarrow[t\to\infty]{}\,0.
  \end{equation*}
  To estimate the second term, we shall exploit the fact that a
  self-similar growth-fragmentation can be constructed from
  a homogeneous one by performing an appropriate Lamperti
  time-substitution on each cell in the system (see
  \cite[Corollary~2]{Bertoin17} or
  \cite[Section~2.1]{Stephenson16}). Specifically, there exists a
  cell system $\CS[Z]\defeq((\CS[Z]_u,\beta_u)\colon u\in\HUT)$, with
  same cumulant function~$\cf$, such that every element in~$\bs X(t)$
  with label $v\in\HUT$ equals $\CS[Z]_u(\tau-\beta_u)$ for some
  $u\in\HUT$ and $\tau\ge0$ fulfilling
  \begin{equation}
    \tau=\int_0^t\biggl(\CS_{\bar v(s)}\!\left(s-b_{\bar v(s)}\right)
      \biggr)^\alpha\,\dd s,%
    \label{eq:A_Galilei}
  \end{equation}
  where $\bar v(s)$ labels the cell in~$\CS$ corresponding to the
  unique ancestor of~$v$ that is alive at time~$s$.
  Further, the connection with compensated fragmentations
  \cite[Proposition~3]{Bertoin17} entails that for every $q\ge0$ with
  $\cf(q)<\infty$,
  \begin{equation*}
    \Ex[1]\!\left[\sum_{u\in\HUT}
    \biggl(\CS[Z]_u\!\left(\tau-\beta_u\right)\biggr)^q\right]
    =\,\exp\bigl(\tau\cf(q)\bigr).
  \end{equation*}
  On the one hand, if~$Z_1(\tau)$ denotes the size of the largest cell
  at time~$\tau$ in~$\CS[Z]$, then Markov's inequality yields
  \begin{equation*}
    \Px[1]\!\left(Z_1(\tau)>\eps^\gamma\right)
      \le\,\eps^{-\gamma q}\,\exp\bigl(\tau\cf(q)\bigr).
  \end{equation*}
  On the other hand, if we purposely take $v\defeq u^*(t)$ then, on the
  complementary event of $\mathcal H_t(\eps)$, we have
  $\CS_{\bar v(s)}(s-b_{\bar v(s)})\ge\eps$ for all
  $s\in\co{0}{t}$ and thus, by~\eqref{eq:A_Galilei},
  $X_1(t)=\CS[Z]_u(\tau-b_u)$ with
  $\tau\ge t\,\eps^\alpha=t^\delta$.
  Hence, fixing $q\in\oo{\cffr}{\cfsr}$ (so that $\cf(q)<0$),
  \begin{equation*}
    \Px[1]\Bigl(\mathcal H_t(\eps)^\complement,\,
      X_1(t)>\eps^\gamma\Bigr)\,
      \le\,\Px[1]\bigl(Z_1(\tau)>\eps^\gamma\bigr)\,
      \le\,t^{q(1-\eta)/\alpha}\,\exp\!\left(t^\delta\cf(q)\right)
      \xrightarrow[t\to\infty]{}\,0.
  \end{equation*}
  Putting the two pieces together we have just showed that, for every
  $0<\eta<1$,
  \begin{equation*}
    \Px[1]\!\left(\frac{\log X_1(t)}{\log t}
      +\frac1\alpha>\frac\eta\alpha\right)
    \xrightarrow[t\to\infty]{}\,0,
  \end{equation*}
  which is the upper bound we wanted.
\end{proof}

\subsection{Freezing the fragmentation}\label{sec:A_Peano}
Suppose now that we ``freeze'' every cell as soon as its size falls
under a fixed diameter $\eps>0$ (which may occur at birth), in
the sense that frozen cells no longer grow or split. To put things
more formally we need a more chronological point of view in the cells
genealogy. For this reason we suppose that the growth-fragmentation has
been constructed as in \cite[Section~2.1]{Stephenson16}, where cells
are now labeled on the infinite binary tree
\begin{equation*}
  \BT\defeq\bigcup_{n=0}^\infty\{1,2\}^n\subset\HUT.
\end{equation*}
Roughly speaking, any jump from a size $x>0$ to some smaller size
$x-y\in\oo{0}{x}$ of a cell with label, say, $u\in\BT$, causes the
death of that cell while at the same time two independent cells labeled
by~$u1$ and~$u2$ are born with initial sizes~$x-y$ and~$y$
respectively. We implicitly reuse the notations of \cref{sec:A_Cartan} within
this new description, e.g.~$\Px$ is the distribution of the cell system
$\CS\defeq(\CS_u\colon u\in\BT)$ when the mother cell starts at size
$x>0$ (i.e.\ has the law~$P_x$). Analogously, $\ell\in\partial\BT$ refers
to a leaf of~$\BT$, and~$\ell[t]$ and~$\CS_\ell(t)$ respectively denote
the label and the process of the unique cell in the branch from~$\es$
to~$\ell$ that is alive at time~$t$. Let us then introduce the first
passage times
\begin{equation*}
  \ft[v]\defeq\inf\,\bigl\{t\ge0\colon\CS_v(t)<\eps\bigr\},
  \qquad v\in\BT\cup\partial\BT,
\end{equation*}
so that the family of frozen cells can be defined as
\begin{equation*}
  \bigl\{\fc{\eps}\bigr\}_{i=1}^\infty
    \defeq\Bigl(\CS_u\bigl(\ft[u]\bigr)\colon u\in\bt\Bigr),
\end{equation*}
with
$\bt\defeq\{u\in\BT\colon u=\ell[\ft[\ell]]\text{ for some }
  \ell\in\partial\BT\}$.
Note that this procedure of freezing cells does not depend on the
scaling parameter~$\alpha$ of the growth-fragmentation
(changing~$\alpha$ just affects the speed at which particles get
frozen). It is proved \cite[Proposition~2.5]{Curien18} that for each
$x>0$, the process of the sum of the sizes of frozen cells raised to
the power~$\cffr$,
\begin{equation*}
  \FM\defeq\!\sum_{i=1}^\infty\fc{\eps}^\cffr,\qquad
  0<\eps\le x,
\end{equation*}
is a backward martingale converging to $M^-(\infty)$ as
$\eps\to0^+$, almost surely and in $\Leb^1(\Px)$.
In the same vein as in~\cite{Martinez05}, we investigate the empirical
measure~$\varphi(\eps)$ defined by
\begin{equation*}
  \bigl\langle\varphi(\eps),f\bigr\rangle\defeq\sum_{i=1}^\infty\fc{\eps}^\cffr\,
  f\!\left(\frac{\fc{\eps}}{\eps}\right)\!.
\end{equation*}
Again, we let~$\tc$ denote the L\'{e}vy process with Laplace
exponent~$\Phi^-$ associated with the pssMp~$\TC$ via Lamperti's
transformation; see~\eqref{eq:A_Banach}. We can check that its L\'{e}vy
measure~$\Lambda^-$ is given by
\begin{equation*}
  \int_\RR g(y)\,\Lambda^-(\dd y)=\int_\RR\left[e^{y\cffr}g(y)
    +\II_{\{y<0\}}(1-e^y)^\cffr g\bigl(\log(1-e^y)\bigr)\right]
  \Lambda(\dd y),
\end{equation*}
see \cite[Theorem~3.9]{Kyprianou14}.
\begin{theorem}
  Suppose~\eqref{eq:A_Liouville}, $\Lambda(\oo{0}{\infty})=0$, and
  that~$\tc$ is not arithmetic.
  Then as $\eps\to0^+$,
  the random measure~$\varphi(\eps)$ converges in $\Px[1]$-probability to
  $M^-(\infty)\,\varphi$, where~$\varphi$ is a deterministic
  probability measure on~$\oo{0}{1}$ specified by
\begin{equation}
  \bigl\langle\varphi,f\bigr\rangle
    \defeq\frac{\cfsr-\cffr}{-\cf'(\cffr)}\iint_{\oo{-\infty}{0}^2}f(e^x)\,
      e^{(\cfsr-\cffr)y}\,\Lambda^-\bigl(\oo{-\infty}{x+y}\bigr)\,%
        \dd x\dd y.%
  \label{eq:A_Maxwell}
\end{equation}
\end{theorem}
\begin{proof}[Proof]
  As said previously we may suppose $\alpha=0$, so that~$\TC$ is just
  the exponential of~$\tc$. After Jagers~\cite{Jagers89}, we can see
  that the random set $\bt\subset\BT$ is a so called
  \emph{optional line} for which the strong branching property holds ---
  intuitively, freezing the cells below~$\eps$ is equivalent to
  freezing those which would descend from a family of cells that have
  first been frozen below $\eps+\delta$, with $\delta>0$ fixed.
  Specifically, by choosing
  $\delta\defeq\sqrt{\eps}-\eps$ for $0<\eps<1$
  and scaling, we can write
\begin{equation*}
  \bigl\langle\varphi(\eps),f\bigr\rangle
    =\sum_{i=1}^\infty\,\underbrace{\vphantom{\sum_j}
       \fc{\sqrt\eps}^\cffr}_{\lambda_i(\eps)}\;
      \underbrace{\sum_{j=1}^\infty\fc[i,j]{\eps_i}^\cffr
        f\!\left(\frac{\fc{\sqrt\eps}\,
          \fc[i,j]{\eps_i}}{\eps}
        \right)\!}_{Y_i(\eps)}\,,
\end{equation*}
  where conditionally on
  $\bs\lambda(\eps)\defeq(\lambda_i(\eps))_{i\ge1}$, the
  $\{\fc[i,j]{\eps_i}\}_{j=1}^\infty,\,i=1,2,\ldots$,
  are independent cell families respectively frozen below
  $\eps_i\defeq\eps/\fc{\sqrt\eps}$. For every
  $1<p<\cfsr/\cffr$, (conditional) Jensen's inequality easily shows that
  the closed martingale~$\FM$ is bounded (by $\Ex[1][M^-(\infty)^p]$)
  in $\Leb^p(\Px[1])$. Hence
\begin{equation*}
  \sup_{0<\eps<1}\Ex[1]\!\left[\left(\sum_{i=1}^\infty
    \lambda_i(\eps)\right)^{\!p\,}\right]<\,\infty,
\end{equation*}
  and, because $\Ex[1][\FM]=\Ex[1][M^-(\infty)]=1$,
\begin{equation*}
  \Ex[1]\!\left[\sum_{i=1}^\infty\lambda_i^p(\eps)\right]
    \le\,\eps^{(p-1)\cffr}\xrightarrow[\eps\to0^+]{}\,0.
\end{equation*}
  The proof then continues like that of~\cref{thm:A_Monge}.
  Similarly to the many-to-one formula, \cref{lem:A_Aristotle} below gives
\begin{equation*}
  \Ex[1]\!\left[\bigl\langle\varphi(\eps),f\bigr\rangle\right]
  =\,\Em[1]\!\left[f\!\left(\TC\bigl(\ft\bigr)/\eps\right)\right]%
    \!,
\end{equation*}
  where $\ft=\inf\{t\ge0\colon\TC(t)<\eps\}$. It thus remains to
  find the distributional limit of
  $\TC(\ft)/\eps$ as $\eps\to0^+$. Observe that up to taking the
  inverse exponential, this random variable corresponds to the overshoot
  above $-\log\eps$ of the spectrally positive L\'{e}vy
  process~$-\tc$, which drifts to~$\infty$ a.s.\ (since
  $E[-\tc(1)]=-(\Phi^-)'(0+)=-\cf'(\cffr)\in\oo{0}{\infty}$). By
  a classical result of renewal theory (see e.g.~\cite{Bingham75} or
  \cite[Theorem~5.7]{Kyprianou14}) we have, for every continuous
  function $g\colon\oo{0}{\infty}\to\RR$ with compact support,
\begin{equation}
  E\!\left[g\bigl(-\tc-(-\log\eps)\bigr)\right]
  \xrightarrow[\eps\to0^+]{}\,\frac1\mu\int_{\oo{0}{\infty}^2}g(x)\,
    \Pi(y+\dd x)\dd y,%
  \label{eq:A_Cauchy}
\end{equation}
  with~$\Pi$ and~$\mu$ respectively the jump measure and the
  expectation at time~$1$ of the ascending ladder height process
  associated with~$-\tc$.
  On the one hand, from \cite[Corollary~4.4.4.(iv)]{Doney07} we get
\begin{equation*}
  \mu=\frac{E[-\tc(1)]}{\mathsf{k}^*},
\end{equation*}
  where~$\mathsf{k}^*$ is the killing rate of the ascending ladder
  height process associated with~$\tc$, and equals the right inverse
  at~$0$ of the Laplace exponent~$\Phi^-$ (see for instance
  \cite[Example~6.11]{Kyprianou14}):
\begin{equation*}
  \mathsf{k}^*\,=\,\sup\left\{t\ge0\colon\Phi^-(t)=0\right\}
    =\,\cfsr-\cffr.
\end{equation*}
  On the other hand, we know since the work of Vigon~\cite{Vigon02}
  (see also \cite[Corollary~7.9]{Kyprianou14}) that~$\Pi$ fulfills
\begin{equation*}
  \Pi\bigl(\oo{y}{\infty}\bigr)
    =\int_0^\infty e^{-\mathsf{k}^*x}
      \Lambda^-\bigl(\oo{-\infty}{-x-y}\bigr)\,\dd x,
  \qquad y>0.
\end{equation*}
  An easy computation then enables us to identify the right-hand sides
  of~\eqref{eq:A_Maxwell} and~\eqref{eq:A_Cauchy}
  (with $g(x)\defeq f(e^{-x})$).
\end{proof}

\begin{lemma}\label{lem:A_Aristotle}
  For every $x>0$ and every bounded measurable function
  $f\colon\oo{0}{\infty}\to\RR$,
\begin{equation*}
  \Ex\!\left[\sum_{u\in\bt}\!\!\!\CS_u^\cffr\bigl(\ft[u]\bigr)\,
    f\Bigl(\CS_u\bigl(\ft[u]\bigr)\Bigr)\right]
  =\,x^\cffr\Em\!\left[f\!\left(\TC\bigl(\ft\bigr)\right)\right]\!,
\end{equation*}
  with the usual convention $f(\cem)\defeq0$.
\end{lemma}
\begin{proof}
  We slightly adapt the proof of~\cite[Proposition~4.1]{Curien18}. To
  this end, recall the intrinsic martingale~$\mathcal{M}^-$ evoked in
  \cref{rem:A_Bernoulli} and, in the paragraph following that remark, the
  definition of the randomly tagged branch~$\RL$. It is here convenient
  to write $u\succ\bt$ if $u\in\BT$ stems from a (unique) node in~$\bt$
  that we then call~$\bar u$ (i.e.\ $\bar u\in\bt$ is a prefix of~$u$).
  The (conditional) distribution of~$\RL$ in~\eqref{eq:A_Kepler} gives
\begin{equation*}
\begin{split}
  \Em\!\left[f\!\left(\TC\bigl(\ft\bigr)\right)\!%
      \II_{\{\RL(k+1)\succ\bt\}}\right]
  =\,\Em\!\left[\sum_{|u|=k+1}\!\!\!\II_{\{u\succ\bt\}}
        \II_{\{u\text{ ancestor of }\RL\}}\,
        f\Bigl(\CS_{\bar u}\bigl(\ft[\bar u]\bigr)\Bigr)\right]\hphantom{\!.}\\[.4em]
  =\,\Em\Bigg[\frac1{\mathcal{M}^-(\infty)}\,\lim_{n\to\infty}
      \sum_{\substack{|u|=k+1\\|v|=n}}\!\!\!\II_{\{u\succ\bt\}}
      \CS_{uv}^\cffr(0)\,f\Bigl(\CS_{\bar u}\bigl(\ft[\bar u]\bigr)\Bigr)
        \Biggr].
\end{split}
\end{equation*}
Rewriting the latter in terms of~$\Px$ simplifies out
$\mathcal{M}^-(\infty)$. The branching property at~$u$ and the
martingale property of~$\mathcal{M}^-$ then entail
\begin{align*}
\Em\!\left[f\!\left(\TC\bigl(\ft\bigr)\right)\!\II_{\{\RL(k+1)\succ\bt\}}
  \right]
   &=\,x^{-\cffr}\Ex\!\left[\sum_{|u|=k+1}\!\!\!\II_{\{u\succ\bt\}}
     \CS_u^\cffr(0)\,f\Bigl(\CS_{\bar u}\bigl(\ft[\bar u]\bigr)\Bigr)
  \right]\!.
\intertext{If we now gather the nodes $u$ which have the same
  ancestor $v\defeq\bar u\in\bt$ and repeat the previous argument, we
  obtain}
\Em\!\left[f\!\left(\TC\bigl(\ft\bigr)\right)\!\II_{\{\RL(k+1)\succ\bt\}}
  \right]
    &=\,x^{-\cffr}\Ex\!\left[\sum_{|v|\le k}\II_{\{v\in\bt\}}
      \CS_v^\cffr\bigl(\ft[v]\bigr)\,
        f\Bigl(\CS_v\bigl(\ft[v]\bigr)\Bigr)\right]\!.
\end{align*}
Since the event $\{\RL(k+1)\succ\bt\}$ must occur for~$k$ large enough
when $\lim_{t\to\infty}\TC(t)=0$ and $f(\cem)=0$ anyway, letting
$k\to\infty$ yields the result by dominated convergence.
\end{proof}

\providecommand{\bysame}{\leavevmode\hbox to3em{\hrulefill}\thinspace}
\providecommand{\MR}{\relax\ifhmode\unskip\space\fi MR }
\providecommand{\MRhref}[2]{%
  \href{http://www.ams.org/mathscinet-getitem?mr=#1}{#2}
}
\providecommand{\href}[2]{#2}
\providecommand{\ARXIV}[1]{\href{http://arXiv.org/abs/#1}{arXiv:#1}}
\providecommand\HAL[2][hal]{\href{http://#1.archives-ouvertes.fr/#1-#2/}{HAL-#2}}

\end{document}